\documentclass[11pt]{article}
\usepackage[margin=1in]{geometry}
\usepackage{amsmath,amssymb,amsthm,pifont,verbatim}
\usepackage{graphicx,subfigure}
\usepackage{epic,fancybox}
\usepackage{epsfig}
\usepackage{epstopdf}
\usepackage{url}
\usepackage{color}
\usepackage{algpseudocode}
\usepackage{algorithm}
\usepackage{scrextend}
\usepackage{hyperref}
\usepackage{tabu}

\newcommand{\spar}{s_\textup{seq}}
\newcommand{\Iop}{\mathcal I}
\newcommand{\Real}{\mathbb{R}}
\newcommand{\Complex}{\mathbb{C}}
\newcommand{\order}{{\mathcal O}}

\newcommand{\hstab}{h_\textup{stab}}
\newcommand{\hacc}{h_\textup{acc}}

\title{A comparison of high order explicit Runge-Kutta, extrapolation, and deferred correction methods in serial and parallel}
\author{David I. Ketcheson\thanks{King Abdullah University of Science \& Technology (KAUST), Thuwal, 23955-6900 Saudi Arabia, \url{david.ketcheson@kaust.edu.sa}} \and Umair bin Waheed\thanks{King Abdullah University of Science \& Technology (KAUST), Thuwal, 23955-6900 Saudi Arabia}}
\begin{document}
\maketitle

\abstract{
We compare the three main types of high-order one-step initial value solvers:
extrapolation, spectral deferred correction, and embedded Runge--Kutta pairs.
We consider orders four through twelve, including both serial and parallel implementations.
We cast extrapolation and deferred correction methods as fixed-order
Runge--Kutta methods, providing a natural framework for the comparison.  
The stability and accuracy properties of the methods are analyzed by theoretical
measures, and these are compared with the results of numerical tests.
In serial, the 8th-order pair of Prince and Dormand (DOP8) is most efficient.
But other high order methods can be more efficient than DOP8 when
implemented in parallel.
This is demonstrated by comparing a parallelized version of the well-known ODEX
code with the (serial) DOP853 code.  For an $N$-body problem with $N=400$, the
experimental extrapolation code is as fast as the tuned Runge--Kutta pair at
loose tolerances, and is up to two times as fast at tight tolerances.

}

\section{Introduction}

The construction of very high order integrators for initial value ordinary differential equations
(ODEs) is challenging:
very high order Runge--Kutta (RK) methods are subject to vast numbers of order conditions,
while very high order linear multistep methods tend to have poor stability properties.
Both extrapolation \cite{Deuflhard1985,hairer1993} and deferred correction
\cite{Daniel_Pereyra_Schumaker_1968,dutt2000} can be used to construct initial
value ODE integrators of arbitrarily high order in a straightforward way.
Both are usually viewed as iterative methods, since they build up
a high order solution based on lower order approximations.  However,
when the order is fixed, methods in both classes can be viewed as Runge--Kutta
methods with a number of stages that grows quadratically with the desired order of accuracy.  

It is natural to ask how these methods compare with standard Runge--Kutta methods.  
Previous studies have compared the relative (serial) efficiency of 
explicit extrapolation and Runge--Kutta (RK) methods \cite{Hull1972,Shampine1986,Hosea1994a},
finding that extrapolation methods have no advantage over moderate to high order
Runge--Kutta methods, and may well be inferior to them
\cite{Shampine1986,Hosea1994a}.  Consequently, extrapolation has recieved
little attention in the last two decades.
It has long been recognized that extrapolation methods offer
excellent opportunities for parallel implementation \cite{Deuflhard1985}.
Nevertheless, to our knowledge no parallel implementation has appeared, and 
comparisons of extrapolation methods have not taken parallel computation into
account, even from a theoretical perspective.
It seems that no work has thoroughly compared the efficiency of spectral deferred correction methods
with that of their extrapolation and RK counterparts.

In this paper we compare the efficiency of
explicit Runge--Kutta, extrapolation, and spectral deferred correction (DC) methods based on their
accuracy and stability properties.  The methods we study are introduced in Section 
\ref{sec:methods} and range in order from four to twelve.
In Section \ref{sec:serial} we give a theoretical analysis based metrics
that are independent of implementation details.
This section is similar in spirit and in methodology to the work of Hosea \& Shampine
\cite{Hosea1994a}.  In Section \ref{sec:test} we validate the theoretical predictions
using simple numerical tests.  
These tests indicate, in agreement with our theoretical analysis and with previous studies,
that extrapolation methods do not have a significant advantage over high order
Runge--Kutta methods, and may in fact be significantly less efficient.
Spectral deferred correction methods generally fare even worse than extrapolation.

In Section \ref{sec:parallel} we analyze the potential of parallel implementations
of extrapolation and deferred correction methods.
We only consider parallelism ``across the method''.
Other approaches to parallelism in time often use parallelism ``across the steps'';
for instance, the parareal algorithm.  Some hybrid approaches include
PFASST \cite{minion2010hybrid,emmett2012toward} and RIDC
\cite{Christlieb_Macdonald_Ong_2010};
see also \cite{guibert2009cyclic}.
Our results should not be used to infer anything about those methods,
since we focus on a simpler approach that does not involve parallelism across multiple steps.

For both extrapolation and (appropriately chosen) deferred correction methods,
the number of stages that must be computed sequentially grows only linearly with
the desired order of accuracy.  Based on simple algorithmic analysis, we extend our 
theoretical analysis to parallel implementations of extrapolation and deferred
correction.  This analysis suggests that extrapolation should be more efficient
than traditional RK methods, at least for computationally intensive problems.
We investigate this further in Section \ref{sec:shared} by performing a 
simple OpenMP parallelization of the ODEX extrapolation code.  The observed
computational speedup is very near the theoretical estimates, and the code
outperforms the DOP853 (serial) code on some test problems.

No study of numerical methods can claim to yield conclusions that are valid for
all possible problems.  Our intent is to give some broadly useful comparisons
and draw general conclusions that can serve as a guide to further studies.
The analysis presented here was performed using the NodePy (Numerical ODEs in Python) package, 
which is freely available from \url{http://github.com/ketch/nodepy}.
Additional code for reproducing experiments in this work can be found at
\url{https://github.com/ketch/high_order_RK_RR}.

\section{High order one-step embedded pairs\label{sec:methods}}
\dictum[P. Deuflhard, 1985]{
...for high order RK formulas the construction of an embedding
RK formula may be beyond human possibilities...}

We are concerned with one-step methods for the solution of the initial value ODE
\begin{align} \label{ivp}
y'(t) & = f(y) & y(t_0) = y_0,
\end{align}
where $y\in\Real^m$, $f: \Real^m \to \Real^m$.
For simplicity of notation, we assume the problem has been written in autonomous form.
An explicit Runge--Kutta pair computes approximations $y_n, \hat{y}_n \approx y(t_n)$
as follows:
\begin{align}
Y_i & = y_n + h \sum_{j=1}^{i-1} a_{ij} f(Y_j) & 1 \le j \le s \\
y_{n+1} & = y_n + h \sum_{j=1}^{s} b_j f(Y_j) \\
\hat{y}_{n+1} & = y_n + h \sum_{j=1}^{s} \hat{b}_j f(Y_j).
\end{align}
Here $h$ is the step size, $s$ denotes number of stages, the {\em stages} $Y_i$ are intermediate approximations, and one evaluation of
$f$ is required for each stage.
The coefficients $A,b,\hat{b}$ determine the accuracy and stability of the method.
The coefficients are typically chosen so that $y_{n+1}$ has local error
$\tau=\order(h^{p})$, and $\hat{y}_{n+1}$ has local error $\hat{\tau}=\order(h^{\hat{p}})$
for some $1<\hat{p}<p$.  Here $p$ is referred to as the order of the method, and sometimes
such a method is referred to as a $p(\hat{p})$ pair.  The value
$\|y_{n+1}-\hat{y}_{n+1}\|$ is used to estimate the error and determine an 
appropriate size for the next step.

The theory of Runge--Kutta order conditions gives
necessary and sufficient conditions for a Runge-Kutta method to be
consistent to a given order \cite{hairer1993,butcher2003}.  For order $p$, these conditions involve
polynomials of degree up to $p$ in the coefficients $A,b$.  The number
of order conditions increases dramatically with $p$: only eight conditions are required
for order four, but order ten requires 1,205 conditions and order fourteen requires 53,263
conditions.  Although the
order conditions possess a great deal of structure and certain simplifying
assumptions can be used to facilitate their solution, the design of efficient Runge--Kutta
pairs of higher than eighth order by direct solution of the order conditions remains
a challenging area.  Some methods of order as high as 14 have been constructed~\cite{feagin2012high}.


\subsection{Extrapolation}
Extrapolation methods provide a straightforward approach to the construction
of high order one-step methods; they can be viewed as Runge--Kutta methods, which
is the approach taken here.
For the mathematical foundations of extrapolation methods we refer the reader to \cite[Section~II.9]{hairer1993}.
The algorithmic structure of extrapolation methods has been considered in detail in previous
works, including \cite{VanderHouwen1990,rauber1997load};
we review the main results here.  Various sequences of step numbers have been proposed, but we consider the
harmonic sequence as it is usually the most efficient \cite{deuflhard_order_1983,Hosea1994a}.
We do not consider the use of smoothing, as previous studies have shown that 
it reduces efficiency \cite{Hosea1994a}.

\subsubsection{Euler extrapolation ({\bf Ex-Euler})}
\begin{figure}
\begin{center}
\includegraphics[height=0.4\textwidth]{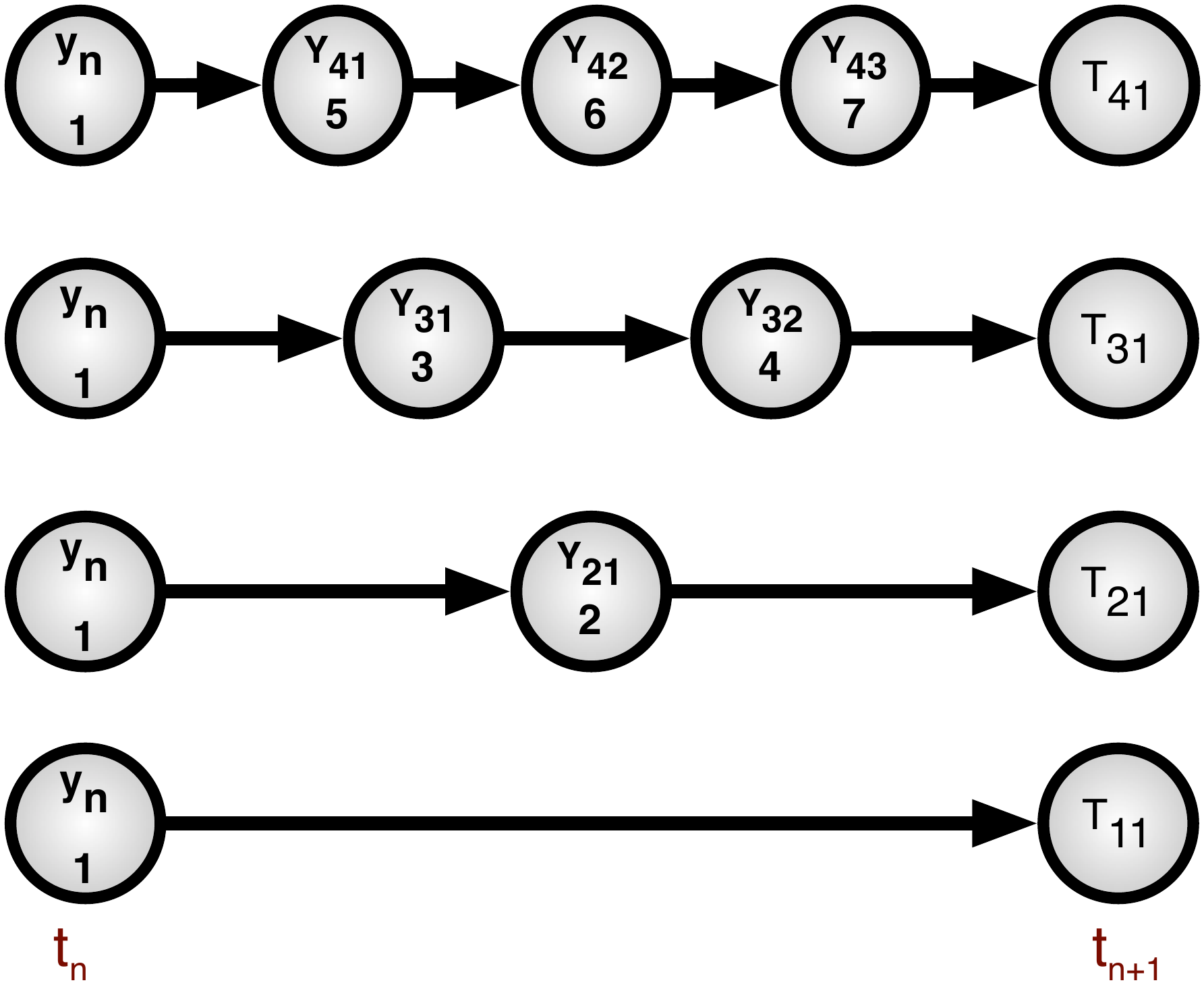}
\caption {Structure of an Euler extrapolation step using the harmonic sequence ${1,2,3,4}$.
Each numbered circle represents a function evaluation,
and the numbers indicate the order in which they are performed.
\label{fig:extrap_structure}}
\end{center}
\end{figure}

Extrapolation is most easily understood by considering the explicit Euler method
\begin{align}
y_{n+1} & = y_n + h f(y_n)
\end{align}
as a building block.
The order $p$ Ex-Euler algorithm computes $p$ approximations to $y(t_{n+1})$
by using the explicit Euler method, first breaking the interval into one step,
then two steps, and so forth.  The approximations to $y(t_{n+1})$ computed in
this manner are all first order accurate and are labeled $T_{11},T_{21}, \dots, T_{p1}$.
These values are combined using the Aitken-Neville interpolation algorithm
to obtain a higher order approximation to $y(t_{n+1})$.  The algorithm is 
depicted in Figure \ref{fig:extrap_structure}.
For error estimation, we use the approximation $T_{p-1,p-1}$ whose accuracy
is one order less.

\begin{algorithm}\caption{Explicit Euler extrapolation ({\bf Ex-Euler})}
\label{alg:extrap}
\begin{algorithmic}

\For{$k = 1 \to p$}  \Comment{Compute first order approximations}
    \State $Y_{k0} = y_n$
    \For{$j=1 \to k$}
        \State $Y_{kj} = Y_{k,j-1} + \frac{h}{k}f(Y_{k,j-1})$
    \EndFor
    \State $T_{k1} = Y_{kk}$
\EndFor

\For{$k=2 \to p$}  \Comment{Extrapolate to get higher order}
    \For{$j=k \to p$}
        \State $T_{jk} = T_{j,k-1} + \frac{T_{j,k-1}-T_{j-1,k-1}}{\frac{j}{j-k+1}-1}$
        \Comment{Aitken-Neville formula for extrapolation to order k}
    \EndFor
\EndFor
\State $y_{n+1} = T_{pp}$ \Comment{New solution value}
\State $\hat{y}_{n+1} = T_{p-1,p-1}$ \Comment{Embedded method solution value}
\end{algorithmic}
\end{algorithm}

Simply counting the number of evaluations of $f$ in Algorithm \ref{alg:extrap}
shows that this is an $s$-stage Runge-Kutta method, where
\begin{align}
s & = \frac{p^2-p+2}{2}.
\end{align}
The quadratic growth of $s$ as the order $p$ is increased leads to relative
inefficiency of very high order extrapolation methods when compared to directly 
constructed Runge--Kutta methods, as we will see in later sections.

\subsubsection{Midpoint extrapolation ({\bf Ex-Midpoint})}
It is common to perform extrapolation based on an integration method whose
error function contains only even terms, such as the midpoint
method~\cite{hairer1993,VanderHouwen1990}.  In this case, each extrapolation
step raises the order of accuracy by two.  We refer to this approach as
Ex-Midpoint and give the algorithm below.  Using midpoint extrapolation to
obtain order $p$ requires about half as many stages, compared to Ex-Euler:
\begin{align}
s & = \frac{p^2+4}{4}.
\end{align}
Again, the number of stages grows quadratically with the order.

\begin{algorithm}\caption{Explicit Midpoint extrapolation ({\bf Ex-Midpoint})}
\label{alg:extrapm}
\begin{algorithmic}
\State $r=p/2$
\For{$k = 1 \to r$}  \Comment{Compute second-order approximations}
    \State $Y_{k0} = y_n$
    \State $Y_{k1} = Y_{k,0} + \frac{h}{2k}f(Y_{k,0})$ \Comment{Initial Euler step}
    \For{$j=2 \to 2k$}
        \State $Y_{kj} = Y_{k,j-2} + \frac{h}{k}f(Y_{k,j-1})$  \Comment{Midpoint steps}
    \EndFor
    \State $T_{k1} = Y_{k,2k}$
\EndFor

\For{$k=2 \to r$}  \Comment{Extrapolate to get higher order}
    \For{$j=k \to r$}
        \State $T_{jk} = T_{j,k-1} + \frac{T_{j,k-1}-T_{j-1,k-1}}{\frac{j^2}{(j-k+1)^2}-1}$
        \Comment{Aitken-Neville formula for extrapolation to order 2k}
    \EndFor
\EndFor
\State $y_{n+1} = T_{rr}$ \Comment{New solution value}
\State $\hat{y}_{n+1} = T_{r-1,r-1}$ \Comment{Embedded method solution value}
\end{algorithmic}
\end{algorithm}

\subsection{Deferred correction ({\bf DC-Euler})}
Like extrapolation, deferred correction has a long history; its application
to initial value problems goes back to \cite{Daniel_Pereyra_Schumaker_1968}.
Recently it has been revived as an area of research, see \cite{dutt2000,Gustafsson_Kress_2001}
and subsequent works.  Here we focus on the class of methods introduced in~\cite{dutt2000},
with a modification introduced in~\cite{liu2008}.
These spectral DC methods are one-step methods and can be constructed for any order of accuracy.

Spectral DC methods start like extrapolation methods, by using a low-order method to
step over subintervals of the time step; the subintervals can be equally sized,
or Chebyshev nodes can be used.  We consider 
methods based on the explicit Euler method and Chebyshev nodes.
Subsequently, high-order polynomial interpolation
of the computed values is used to approximate the integral of the error, or defect.
Then the method steps over the same nodes again, applying a correction.  This procedure
is repeated until the desired accuracy is achieved.  

\begin{figure}
\begin{center}
\includegraphics[height=0.4\textwidth]{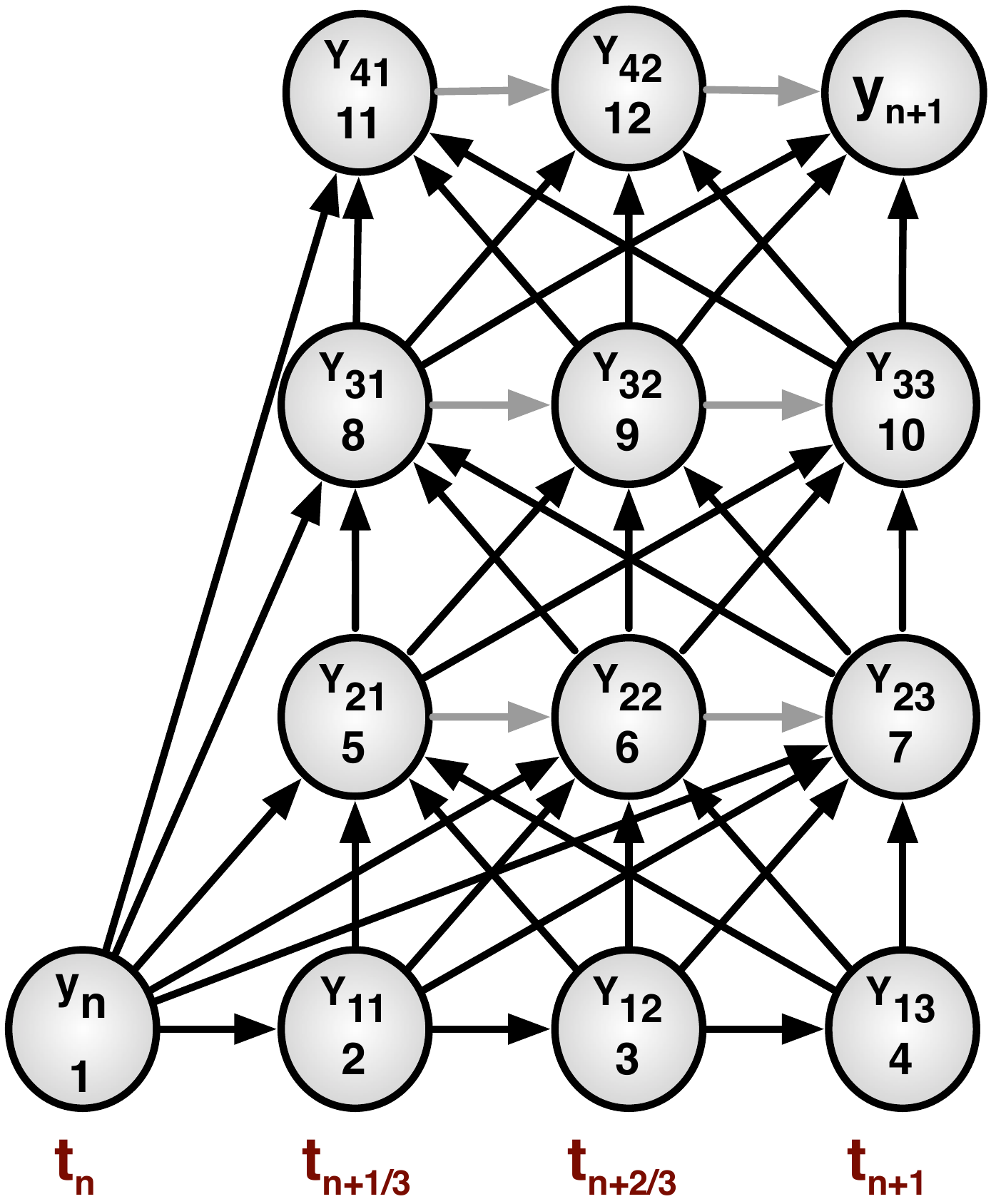}
\caption {Structure of a 4th-order spectral DC step using 3 Euler substeps.
Each numbered circle represents a function evaluation,
and the numbers indicate the order in which they are performed.
The black arrows represent dependencies; the grey arrows are dependencies
that vanish when $\theta=0$.  Note that node 1 is connected to all
other nodes; some of those arrows have been omitted for clarity.
Thus the solution at each node depends on all solutions from the previous
iteration and, unless $\theta=0$, on its predecessor in the current iteration.
\label{fig:idc_structure}}
\end{center}
\end{figure}

A modification of the spectral DC method appears in \cite{liu2008}, in which a parameter $\theta$
is used to adjust the dependence of the correction steps on previous iterations.  
The original scheme corresponds to $\theta=1$; by taking $\theta\in[0,1]$ the
stability of the method can be improved.  
Given a fixed order of accuracy and a predictor method, the resulting spectral DC method can be
written as a Runge--Kutta method \cite{gottlieb2009}.
The algorithm is defined below (the values $c_j$ denote the locations of the Chebyshev
nodes) and depicted in Figure \ref{fig:idc_structure}.
For error estimation, we use the solution from the next-to-last correction
iteration, whose order is one less than that of the overall method.

\begin{algorithm}\caption{Explicit Euler-based deferred correction ({\bf DC-Euler})}
\label{alg:idc}
\begin{algorithmic}

\State $Y_{10} = y_n$
\For{$k = 1 \to p-1$} \Comment{Compute initial prediction}
    \State $Y_{1k} = Y_{1,k-1} + (c_{k+1}-c_k)h f(Y_{1,k-1})$
\EndFor

\For{$k=2 \to p$}  \Comment{Compute successive corrections}
    \State $Y_{k0} = y_n$
    \For{$j=1 \to p-1$}
        \State $Y_{kj} = Y_{k,j-1} + h\theta (f(Y_{k,j-1})-f(Y_{k-1,j-1})) + \Iop_{j-1}^j(f(Y_{k-1,:}))$
    \EndFor
\EndFor
\State $y_{n+1} = Y_{p,p-1}$ \Comment{New solution value}
\State $\hat{y}_{n+1} = Y_{p-1,p-1}$ \Comment{New solution value}
\end{algorithmic}
\end{algorithm}

In Algorithm \ref{alg:idc}, $\Iop_{j-1}^j(f(Y_{k-1,:}))$ represents the integral of
the degree $p-1$ polynomial that interpolates the points $Y_{k-1,j}$ for $j=1,\dots,p-1$,
over the interval $[t_n + c_j h, t_n + c_{j+1} h]$.  Thus, for $\theta=0$, the algorithm
becomes a discrete version of Picard iteration.

The number of stages per step is
\begin{align}
    s = p(p-1)
\end{align}
unless $\theta=0$, in which case the stages $Y_{p,j}$ (for $j<p-1$) need not be computed at
all since $Y_{p,p-1}$ depends only on the $Y_{p-1,j}$.  Then the number of stages
per step reduces to $(p-1)^2+1$.


\subsection{Reference Runge--Kutta methods}
In this work we use the following existing Runge-Kutta pairs as benchmarks for evaluating
extrapolation and deferred correction methods:
\begin{itemize}
\item Fourth order: the embedded formula of Merson 4(3)~\cite[pg. 167]{hairer1993}
\item Sixth order: the 6(5) pair of Calvo et. al.~\cite{calvo_new_1990}, which was
found to be the most efficient out of those considered by Hosea and Shampine~\cite{Hosea1994a}
\item Eighth order: the well-known Prince-Dormand 8(7) pair~\cite{Prince1981}
\item Tenth order: the 10(8) pair of Curtis~\cite{curtis_high-order_1975}
\item Twelfth order: The 12(9) pair of Ono~\cite{ono200625}
\end{itemize}
It should be stressed that finding pairs of order higher than eight
is still very challenging, and the tenth- and twelfth-order pairs here are not 
expected to be as efficient as that of Prince-Dormand.

\section{Concurrency\label{sec:parallel}}
\dictum[P. Deuflhard, 1985]{
In view of an implementation on parallel computers, extrapolation methods (as
opposed to RKp methods or multistep methods) have an important distinguishing
feature: the rows can be computed independently.}

If a Runge-Kutta method includes stages that are mutually independent, then those stages
may be computed concurrently \cite{jackson}.
In this section we investigate theoretically achievable parallel speedup
and efficiency of extrapolation and deferred correction methods.
Our goal is to determine hardware- and problem- independent upper bounds based purely on
algorithmic concerns.  We do not attempt to account for machine-specific overhead or communication,
although the simple parallel tests in Section \ref{sec:odex} suggest that the bounds
we give are realistically achievable for at least some classes of moderate-sized problems.
Previous works that have considered concurrency in explicit extrapolation and deferred correction methods
include \cite{simonsen1990extrapolation,VanderHouwen1990,rauber1997load,burrage1995parallel,guibert2009cyclic,emmett2012toward,Kappeller_Kiehl_Perzl_Lenke_1996,lustman1992solution,minion2010hybrid}.

 \subsection{Computational model and speedup\label{sec:speedup}}
    As in the serial case, our computational model is based on the assumption that evaluation
    of $f$ is sufficiently expensive so that all other operations (e.g., arithmetic,
    step size selection) are negligible by comparison.

\begin{figure}
    \begin{center}
    \includegraphics[height=0.4\textwidth]{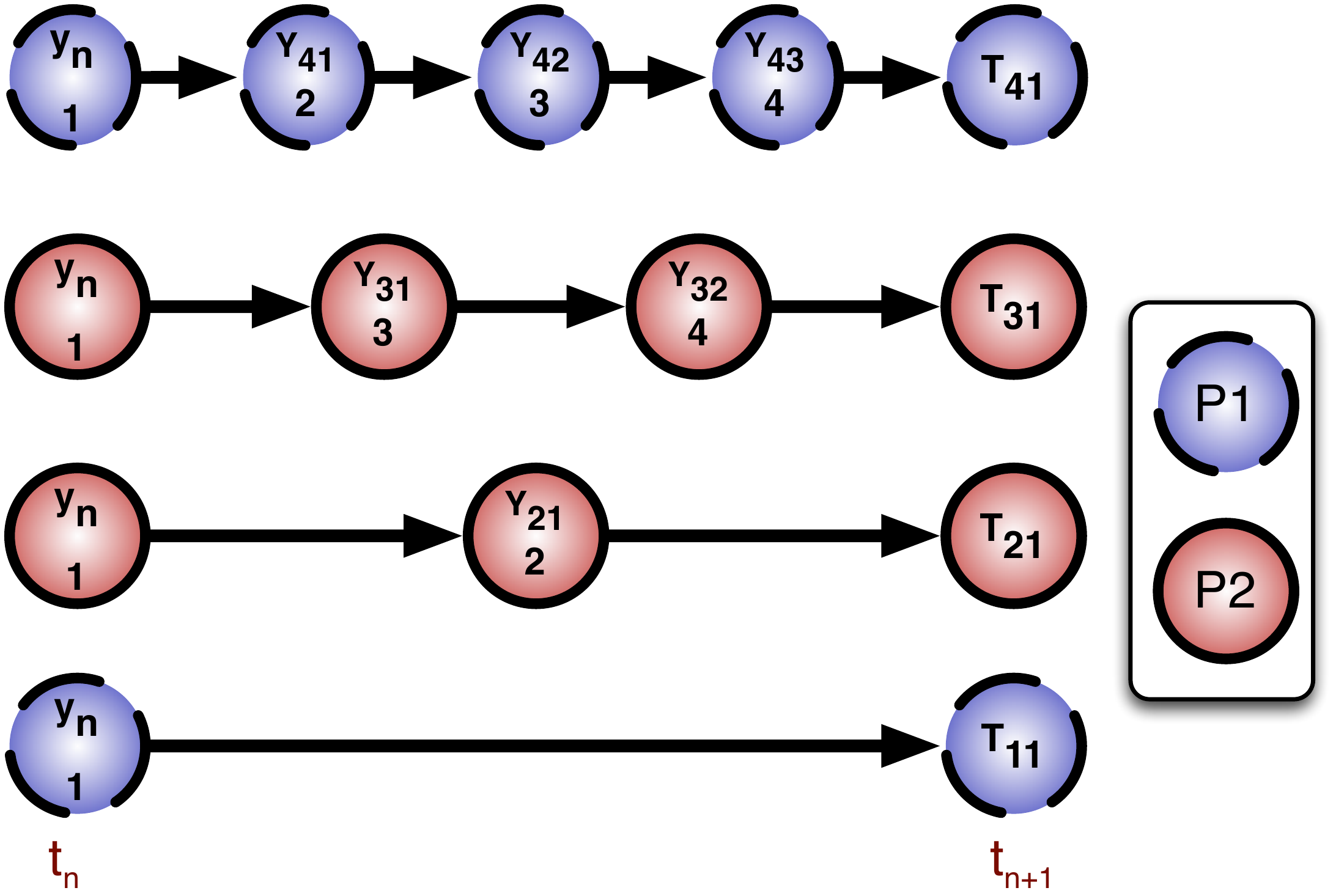}
    \caption {Exploiting concurrency in an Euler extrapolation step using 2 processes.
    The blue circles with broken border are computed by process 1 and the red
    circles with solid border are computed by process 2.  Observe that only $\spar=4$
    sequential function evaluations are required for each process, as opposed to
    the $s=7$ sequential evaluations required in serial.
    \label{fig:extrap_2proc}}
    \end{center}
    \end{figure}

    Typically, stage $y_j$ of an explicit Runge--Kutta method depends on all the 
    previous stages $y_1, y_2, \dots, y_{j-1}$.  
    However, if $y_j$ does not depend on $y_{j-1}$, then these two stages may be computed
    simultaneously on a parallel computer.  More generally, by interpreting the incidence matrix of $A$ as the 
    adjacency matrix of a directed graph $G(A)$, one can determine precisely which stages
    may be computed concurrently and how much speedup may be achieved.
    For extrapolation methods, the computation of each $T_{k1}$ may be performed independently
    in parallel \cite{Deuflhard1985}, as depicted in Figure \ref{fig:extrap_2proc}.  
    Unlike some previous authors, we do not consider parallel implementation of
    the extrapolation process (i.e., the second loop in Algorithm 1) since it
    does not include any evaluations of $f$ (so our computational model assumes its
    cost is negligible anyway).

\begin{figure}
    \begin{center}
    \includegraphics[height=0.4\textwidth]{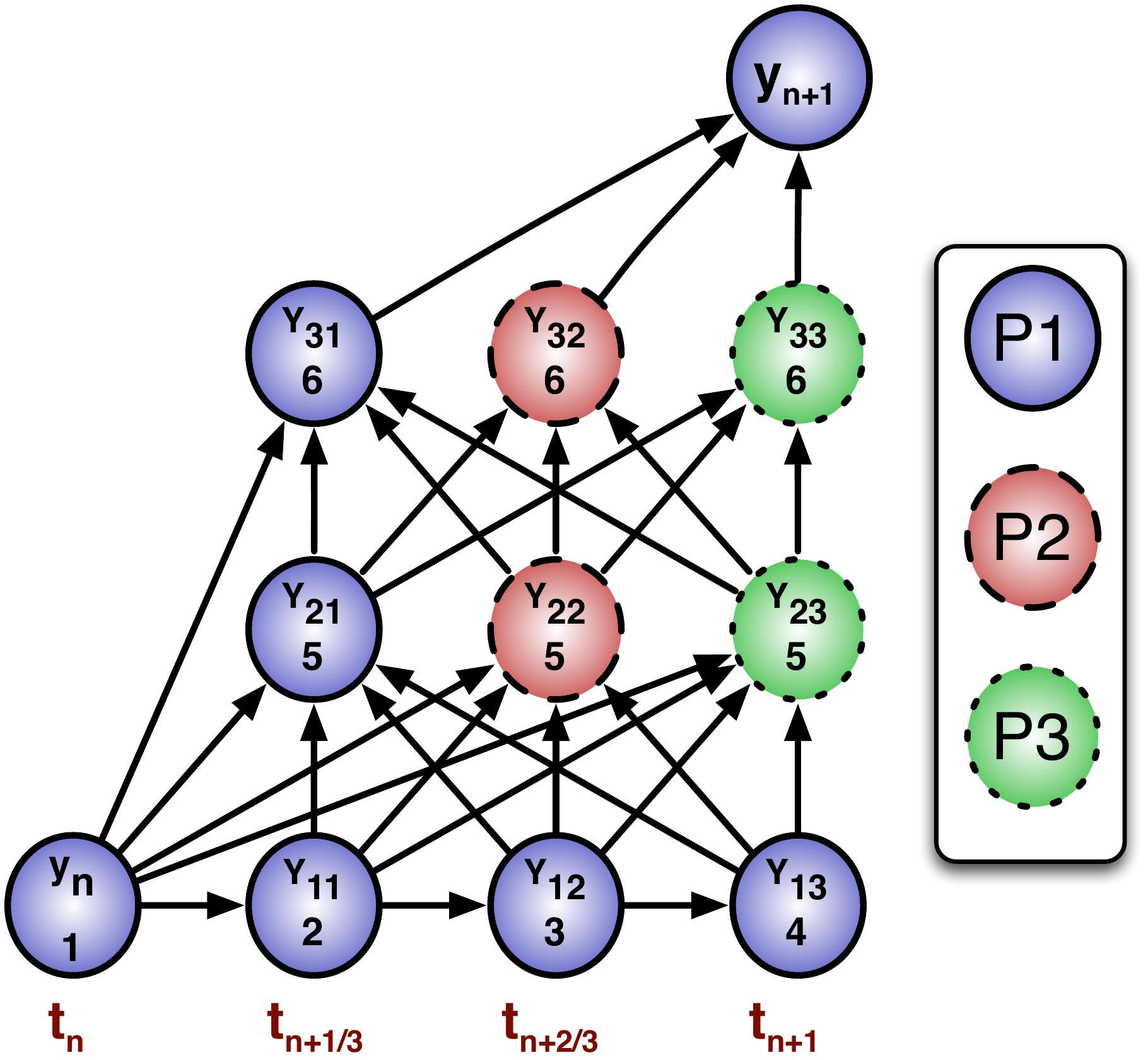}
    \caption {Exploiting concurrency in a 4th-order spectral DC step (with $\theta=0$)
    using 3 Euler substeps and 3 processes.
    The color and border of each circle indicate which process evaluates it.
    Observe that only $\spar=6$
    sequential function evaluations are required for each process, as opposed to
    the $s=10$ sequential evaluations required in serial (12 in serial when $\theta\ne 0$).
    Note that node 1 is connected to all
    other nodes; some of those arrows have been omitted for clarity.
    More synchronization is required than for a similar extrapolation step.
    \label{fig:dc_3proc}}
    \end{center}
    \end{figure}

    For the deferred correction methods we consider, parallel computation is advantageous
    only if $\theta=0$; the resulting parallel algorithm is depicted in
    Figure~\ref{fig:dc_3proc}.  A different approach
    to parallelism in DC methods is taken by the RIDC method \cite{Christlieb_Macdonald_Ong_2010};
    see also \cite{guibert2009cyclic}.
    Deferred correction has also been combined with the parareal algorithm to achieve
    parallel speedup \cite{minion2010hybrid,emmett2012toward}.

    We define the {\em minimum number of sequential stages} $\spar$ as the minimum
    number of sequential function evaluations that must be made when parallelism
    is taken into account.  To make this more precise, let us label each node in 
    the graph $G(A)$ by the index of the stage it corresponds to, with the node
    corresponding to $y_{n+1}$ labeled $s+1$.  Then
    \begin{align}
        \spar = \max_j \{  \text{ path length from node $1$ to node $s+1$}\}.
    \end{align}
    The quantity $\spar$ represents the minimum time required to
    take one step with a given method on a parallel computer, 
    in units of the cost of a single derivative evaluation.
    For instance, the maximum path length for the method shown in Figure~\ref{fig:extrap_2proc} is equal to 4;
    for the method in Figure~\ref{fig:dc_3proc} it is 6.
    The maximum potential parallel speedup is 
    \begin{align}
    S & = s/\spar.
    \end{align}
    The minimum number of processes required to achieve speedup $S$ is denoted by $P$
    (equivalently, $P$ is the maximum number of processes that
    can usefully be employed by the method).  Finally, let $E$ denote the theoretical parallel
    efficiency (here we use the term in the sense that is common in the parallel computing literature)
    that could be achieved by spreading the computation over $P$ processes:
    \begin{align}
        E = \frac{s}{P\spar} = \frac{S}{P}.
    \end{align}
    Note that $E$ is an upper bound on the achievable parallel efficiency; it accounts only
    for inefficiencies due to load imbalancing.
    It does not, of course, account for additional
    implementation-dependent losses in efficiency due to overhead or communication.

    Table \ref{tbl:par} shows the parallel algorithmic properties of fixed-order extrapolation
    and deferred correction methods.   
    Note that for deferred correction methods with $\theta\ne 0$, we have
    $\spar=s$, i.e., no parallel computation of stages is possible.

    \begin{table}
    \centering
    {\tabulinesep=0.5mm
    \begin{tabu}{c|c|c|c|c|c}
    Method & $s$ & $\spar$ & $S$ & $P$ & $E$ \\ \hline
    Ex-Euler & $\frac{p^2-p+2}{2}$ & $p$ & $\frac{p^2-p+2}{2p}$  & $\lceil \frac{p}{2}\rceil$ & $\frac{p^2-p+2}{2p\lceil \frac{p}{2}\rceil}$ \\ \hline
    Ex-Midpoint & $\frac{p^2+4}{4}$ & $p$ & $\frac{p^2+4}{4p}$ & $\lceil \frac{p+2}{4}\rceil$  & $\frac{p^2+4}{4p\lceil \frac{p+2}{4}\rceil}$ \\ \hline
    DC-Euler, $\theta=0$ & $(p-1)^2+1$ & $2(p-1)$ & $\frac{(p-1)^2+1}{2(p-1)}$ & $p-1$ & $\frac{(p-1)^2+1}{2(p-1)}$\\ \hline
    DC-Euler, $\theta \ne 0$ & $p(p-1)$ & $p(p-1)$ & 1 & 1 & - \\
    \end{tabu}}
    \caption{Parallel implementation properties of extrapolation and deferred correction methods.
    $s$: number of stages; $\spar$: number of sequentially dependent stages;
    $S=s/\spar$: optimal speedup; $P$: number of processes required to achieve optimal
    speedup; $E = S/P$: parallel efficiency.
    \label{tbl:par}}
    \end{table}

To our knowledge, no parallel implementation has been made of the 
deferred correction methods we consider here.  However, the parallel iterated
RK methods of \cite{van1990parallel} have a similar flavor.  For parallel implementation
of a {\em revisionist} DC method, see \cite{Christlieb_Macdonald_Ong_2010}.

\section{Theoretical measures of efficiency\label{sec:serial}}

    Here we describe the theoretical metrics we use to evaluate the methods.
    Our metrics are fairly standard; a useful and thorough reference is \cite{Kennedy2000}.
    The overarching metric for comparing methods is efficiency: the number of function evaluations
    required to integrate a given problem over a specified time interval to a specified
    accuracy.  We assume that function evaluations are relatively expensive so that
    other arithmetic operations and overhead for things like step size selection are not significant.

    The number of function evaluations is the product of the number of stages of the method
    and the number of steps that must be taken.  The number of steps to be
    taken depends on the step size $h$, which is usually determined adaptively to
    satisfy accuracy and stability constraints, i.e.
    \begin{align}
        h = \min(\hstab,\hacc)
    \end{align}
    where $\hstab, \hacc$ are the maximum step sizes that ensure numerical stability
    and prescribed accuracy, respectively.  Since the cost of a step is proportional
    to the number of stages of the method, $s$, then a fair measure of efficiency is
    $h/s$.  A simple observation
    that partially explains results in this section is as follows: extrapolation
    and deferred correction are straightforward approaches to creating methods
    that satisfy the huge numbers of order conditions for very high order Runge--Kutta
    methods.  However, this straightforward approach comes with a cost: they
    use many more than the minimum necessary number of stages to achieve a particular order,
    leading to relatively low efficiency.

 \subsection{Absolute stability}
    The stable step size $\hstab$ is typically the limiting factor when a very loose
    error tolerance is applied.
    A method's region of absolute stability (in conjunction with the spectrum of $f'$)
    typically dictates $\hstab$.

    In order to make broad comparisons, we measure the size
    of the and real-axis interval that is contained in the absolute stability
    region.  Specifically, let $S\subset \Complex$ denote the region of absolute
    stability; then we measure
    \begin{align}
        I_\textup{real} & = \max\{r\ge 0 : [-r,0]\subset S\} \\
        I_\textup{imag} & = \max\{r\ge 0 : [-ir,ir]\subset S\}.
    \end{align}

    Determination of the stability region for very high order methods can
    be numerically delicate; for instance, the stability function for the
    8th-order deferred correction method is a polynomial of degree 56!
    Because of this, all stability calculations presented here have been
    performed using exact (rational) arithmetic, not in floating point.

\begin{figure}[ht!]
    \begin{center}
    \subfigure[Real Stability Interval]{%
    \label{fig:stability1}
    \includegraphics[width=0.4\textwidth]{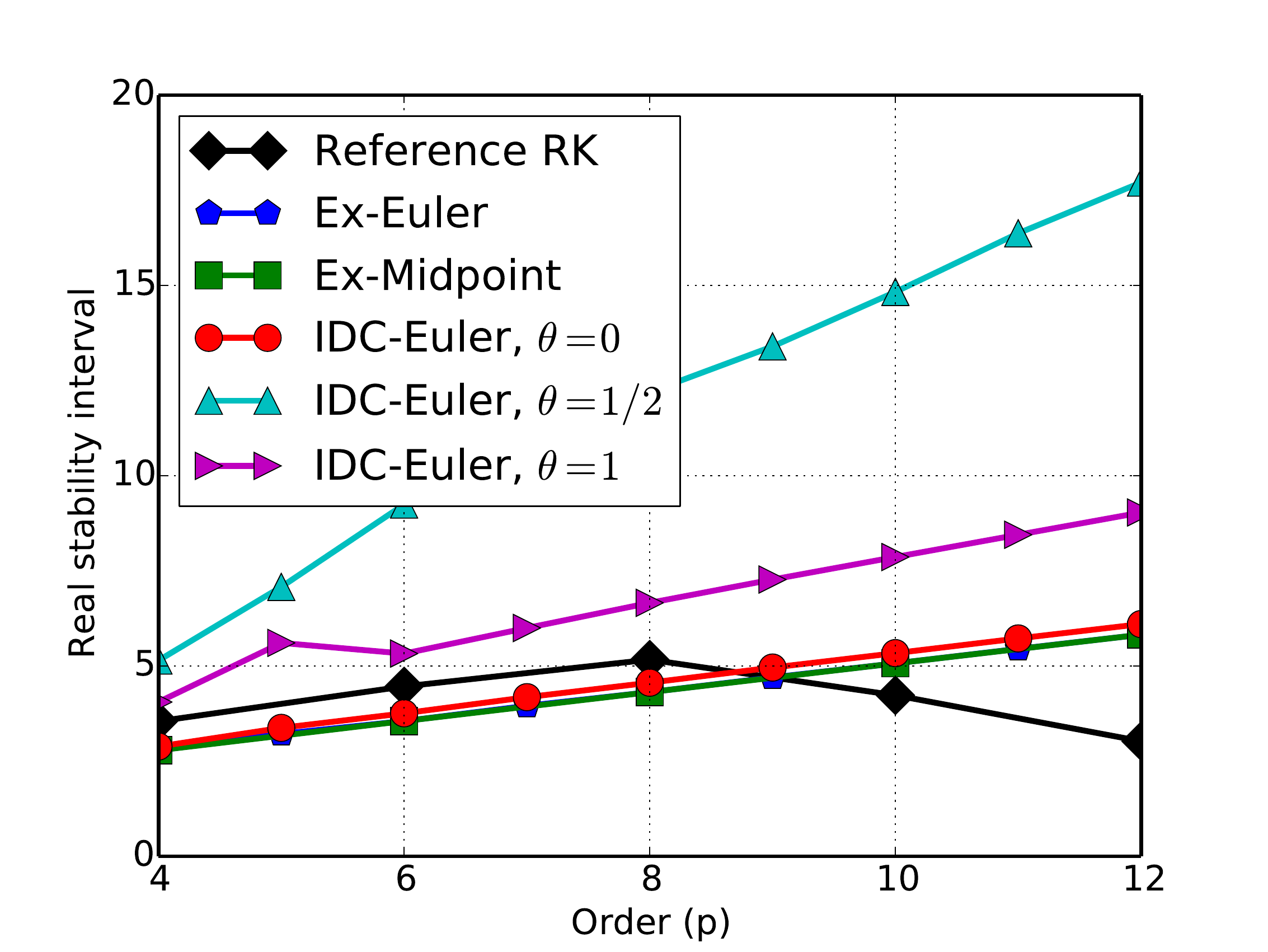}
    }%
    \subfigure[Scaled Real Stability Interval]{%
    \label{fig:scaledstability1}
    \includegraphics[width=0.4\textwidth]{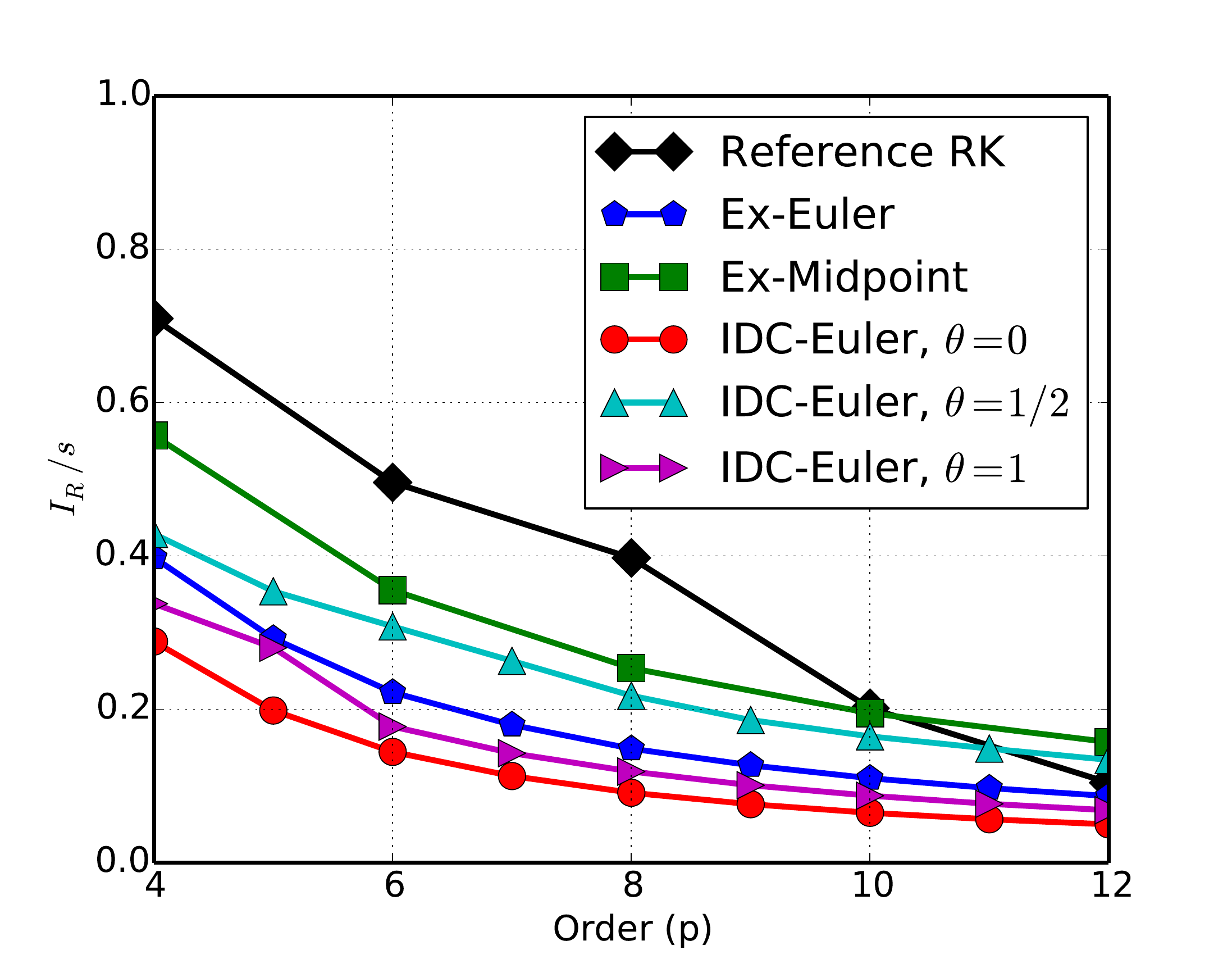}
    }%
    \end{center}
    \caption{%
        Comparison of stability regions for reference methods, Euler
        extrapolation, midpoint extrapolation and deferred correction. 
        \label{fig:stability}}
    \end{figure}

    Figure~\ref{fig:stability1} and Table~\ref{tbl:stability} show real and imaginary
    stability interval sizes 
    for Ex-Euler, Ex-Midpoint, and DC-Euler methods of orders 4-12. 
    We show the real stability intervals of the deferred correction methods
    with three different values of $\theta$, because this interval has a strong
    dependence on $\theta$.
    For all classes of methods, the overall size of the stability region grows
    with increasing order.  However, many methods have $I_\textup{imag}=0$.
    Note that the stability regions for Ex-Euler and Ex-Midpoint are identical
    since both have stability polynomial
    \begin{align}
    \sum_{k=0}^p \frac{z^p}{p!},
    \end{align}
    i.e., the degree-$p$ Taylor polynomial of the exponential function.

    A fair metric for efficiency is obtained by dividing these interval sizes
    by the number of stages in the method.  The result is shown in 
    Figure~\ref{fig:scaledstability1}.  Higher-order methods have smaller
    relative stability regions.
    For orders $p\le 10$, the reference RK methods have
    better real stability properties.
    We caution that, for high order methods, the boundary of the stability region typically
    lies very close to the imaginary axis, so values of the amplification factor may 
    differ from unity by less than roundoff over a large interval.
    For instance, the 10th-order extrapolation method has $I_\textup{imag}=0$,
    but the magnitude of its stability polynomial differs from unity by less than
    $1.4\times 10^{-15}$ over the interval $[-i/4,i/4]$.
    It is not clear whether precise
    measures of $I_\textup{imag}$ are relevant for such methods in practical situations.
    
    Here for simplicity we have considered only the stability region of the principal method;
    in the design of embedded pairs, it is important that the embedded method have a similar
    stability region.  All the pairs considered here seem to have fairly well matched stability
    regions.


    \begin{table}
    \centering
    \begin{tabular}{r|l|l|l|l}
    Order & Reference RK & Ex-Euler & Ex-Midpoint & DC-Euler, $\theta=0$ \\ \hline
    4     & 3.46         & 2.83     & 2.83        & 2.93     \\ \hline
    5     & -            & 0        & -           & 0        \\ \hline
    6     & 2.61         & 0        & 0           & 0        \\ \hline
    7     & -            & 1.76     & -           & 1.82     \\ \hline
    8     & 0            & 3.40     & 3.40        & 3.52     \\ \hline
    9     & -            & 0        & -           & 0        \\ \hline
    10    & 0            & 0        & 0           & 0        \\ \hline
    11    & -            & 1.70     & -           & 1.75     \\ \hline
    \end{tabular}
    \caption{Imaginary stability intervals.\label{tbl:stability}}
    \end{table}


 \subsection{Accuracy efficiency}


    Typically, the local error is controlled
    by requiring that $\|y_{n+1}-\hat{y}_{n+1}\|<\epsilon$ for some tolerance $\epsilon>0$.
    When the maximum stable step size does not yield sufficient accuracy, the accuracy 
    constraint determines the step size.  This is typically the case when
    the error tolerance is reasonably small.  In theoretical analyses, the {\em
    principal error norm} \cite{Kennedy2000} 
    \begin{align} \label{error-norm}
    C_{p+1} & = \left( \sum_k (\tau_k^{(p+1)})^2 \right)^{\frac{1}{2}}
    \end{align}
    is often used as a way to compare accuracy between two methods of the same order.
    Here the constants $\tau_k^{(p+1)}$ are the coefficients appearing in the leading
    order truncation error terms.

    Assuming that the one-step error is proportional to $C_{p+1}h^{p+1}$ leads
    to a fair comparison of accuracy efficiency given by the {\em
    accuracy efficiency index}, introduced in \cite{Hosea1994a}:
    \begin{equation}
    \eta = \frac{1}{s} \left(\frac{1}{C_{p+1}}\right)^{1/p+1}.
    \label{eq:acceff}
    \end{equation}

\begin{figure}
    \begin{center}
    \subfigure[Serial accuracy efficiency\label{fig:acc_eff}]{%
    \includegraphics[width=0.4\textwidth]{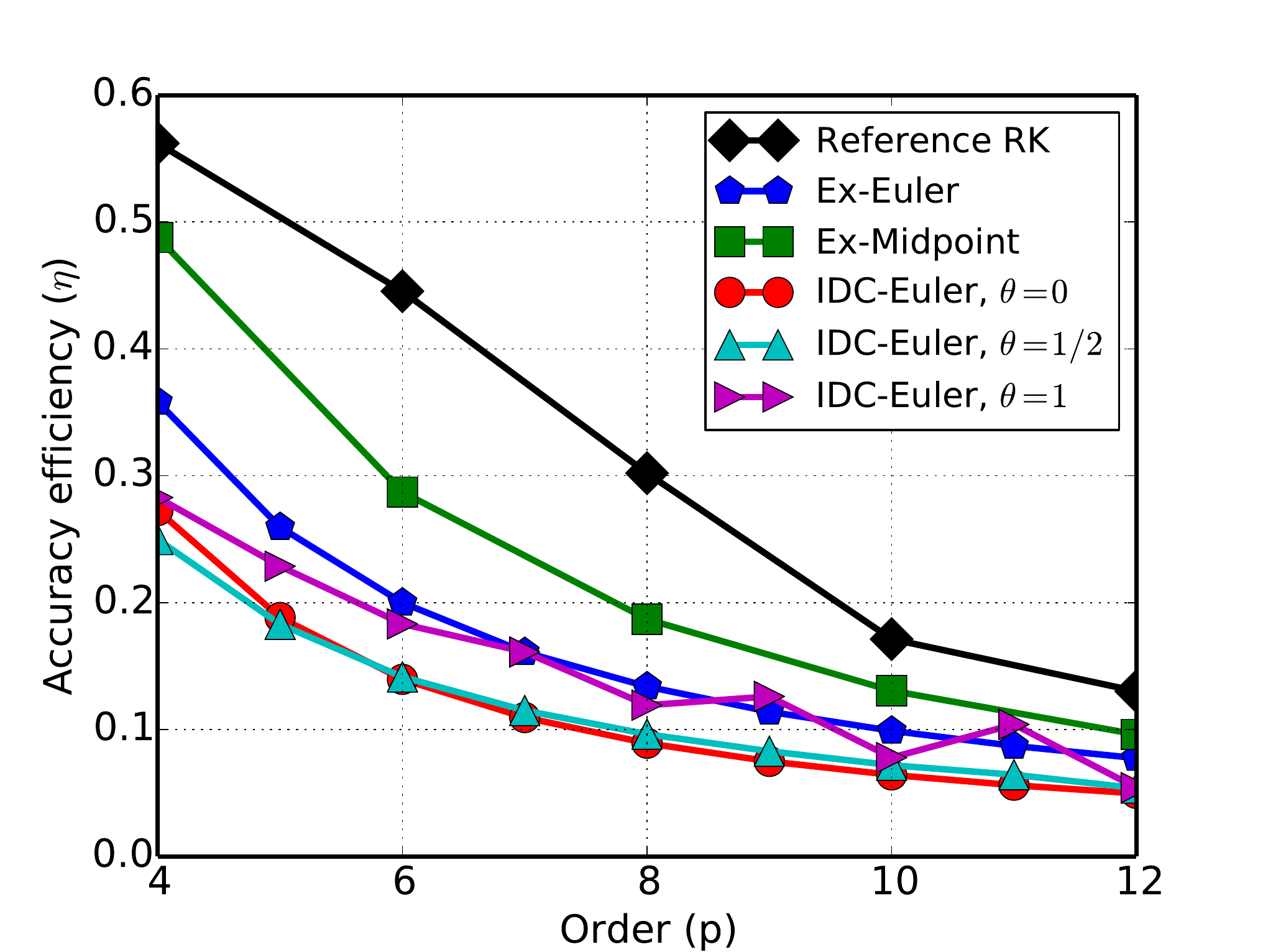}}
    \subfigure[Ideal parallel accuracy efficiency\label{fig:acc_eff_p}]{%
    \includegraphics[width=0.4\textwidth]{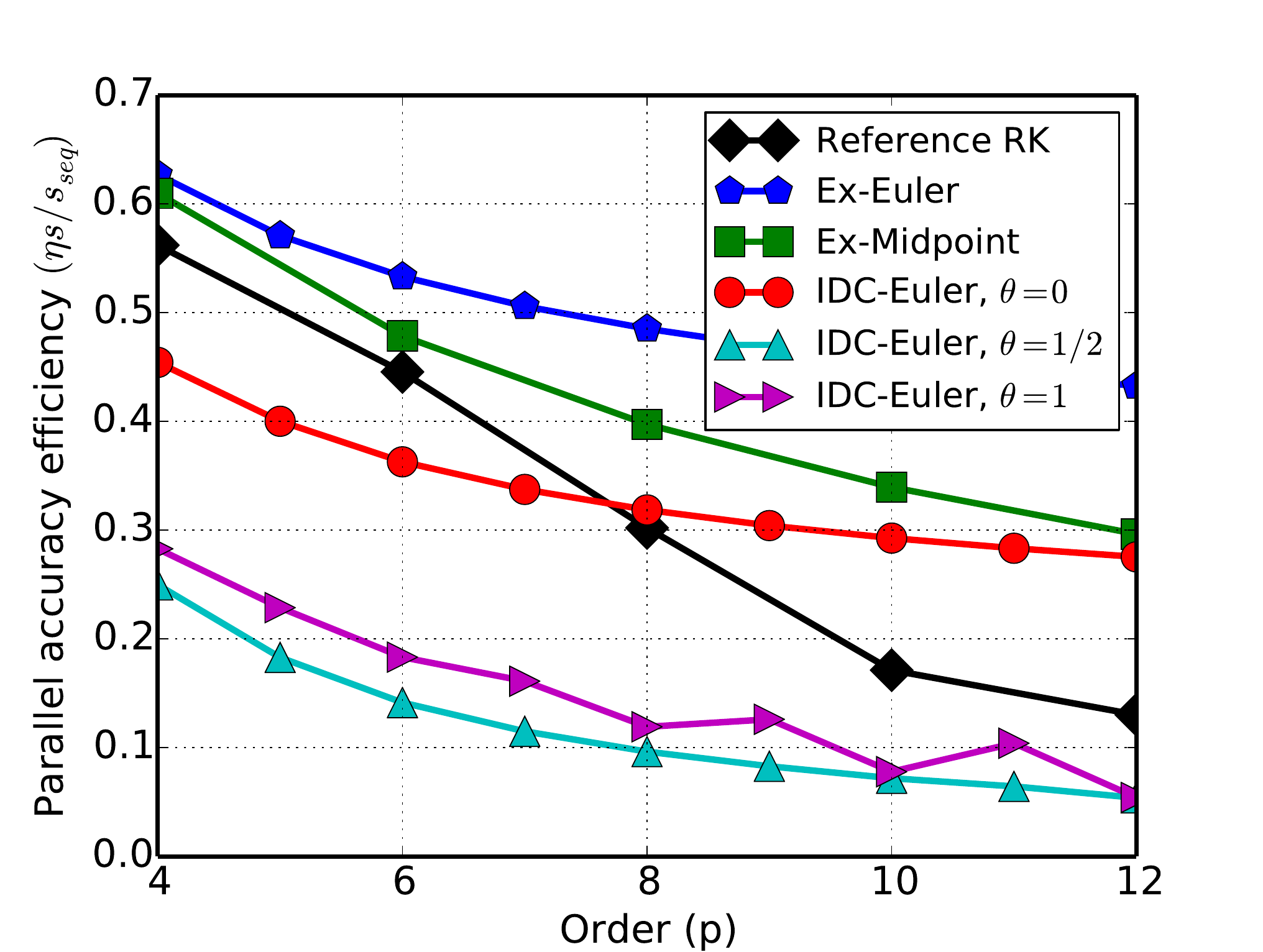}}
    \caption{Accuracy efficiency}
    \end{center}
    \end{figure}

    Figure~\ref{fig:acc_eff} plots the accuracy efficiency index for the methods under consideration.
    Interestingly, a ranking of methods based on this metric gives the same ordering as that based on
    $I_\textup{real}/s$.

\subsection{Accuracy and stability metrics}
In order to determine idealized accuracy and stability efficiency measures, we take
the speedup factor $s/\spar$ into account.  In other words, we consider
\begin{equation}
\frac{s}{\spar} \eta = \frac{1}{\spar} \left(\frac{1}{C_{p+1}}\right)^{1/p+1},
\label{eq:acceffp}
\end{equation}
as a measure of accuracy efficiency.  A similar scaling could be used to study 
stability efficiency of parallel implementations.
We stress that in this context {\em efficiency} relates to the number of function evaluations
required to advance to a given time, and is not related to the usual concept of parallel efficiency.

Figure~\ref{fig:acc_eff_p} shows the accuracy efficiency, rescaled by the speedup factor.
Comparing with Figure~\ref{fig:acc_eff}, we see a very different picture for methods of order
8 and above.  Extrapolation methods are the most efficient, while the reference RK 
methods give the weakest showing -- since they do not benefit from parallelism.

\subsection{Predictions}
    The theoretical measures above indicate that fixed-order extrapolation and deferred correction
    methods are less efficient than traditional Runge--Kutta methods, at least up to
    order eight.  At higher orders, the disadvantage of extrapolation and spectral DC are less
    pronounced, but they still offer no theoretical advantage.  When parallelism is taken
    into account, extrapolation and deferred correction offer a significant theoretical advantage.

\section{Performance tests\label{sec:test}}
In this section we perform numerical tests, solving some initial value problems with
the methods under consideration, to validate the theoretical predictions of the last section.

In addition to the tests shown, we tested all methods on a collection of
problems known as the non-stiff DETEST suite~\cite{Hull1972}.  The results (not
shown here) are broadly consistent with those seen in the test problems below.

\subsection{Verification tests}
For each of the pairs considered, we performed convergence tests using a sequence
of fixed step sizes with several nonlinear systems of ODEs, in order to verify that the
expected rate of convergence is achieved in practice.  We also checked that the
coefficients of each method satisfy the order conditions exactly (in rational arithmetic).

\subsection{Step size control}
For step size selection, we use a standard $I$-controller~\cite{Kennedy2000}:
\begin{equation}
h_{n+1}^* = \kappa h_n \left(\frac{\epsilon}{||\delta_{n+1}||_\infty}\right)^\alpha.
\end{equation}
Here $\epsilon$ is the chosen integration tolerance and
$\delta_{n+1}=y_{n+1}-\hat{y}_{n+1}$. We
take $\kappa=0.9$ and $\alpha=0.7/p$, where $p$ is the order of the embedded method.
The step size is not allowed to increase or decrease too suddenly; we use~\cite{hairer1993}:
\begin{equation}
h_{n+1}= \text{min}\left(\kappa_{max}h_n,\text{max}\left(\kappa_{min}h_n,h^*_{n+1}\right)\right)
\end{equation}
with $\kappa_{min}=0.2$ and $\kappa_{max}=5$.
A step is rejected if the error estimate exceeds the tolerance;
i.e., if $\|\delta_n\|_\infty>\epsilon$.

All tests in this work were also run with a $PI$-controller, and very similar
results were obtained.

\subsection{Test problems and results\label{sec:serial-results}}
\subsubsection{Three-body problem\label{threebody}}
\begin{figure}
\begin{center}
\subfigure[6th order]{\label{fig:d1_6}
\includegraphics[width=0.4\textwidth,height=0.32\textwidth]{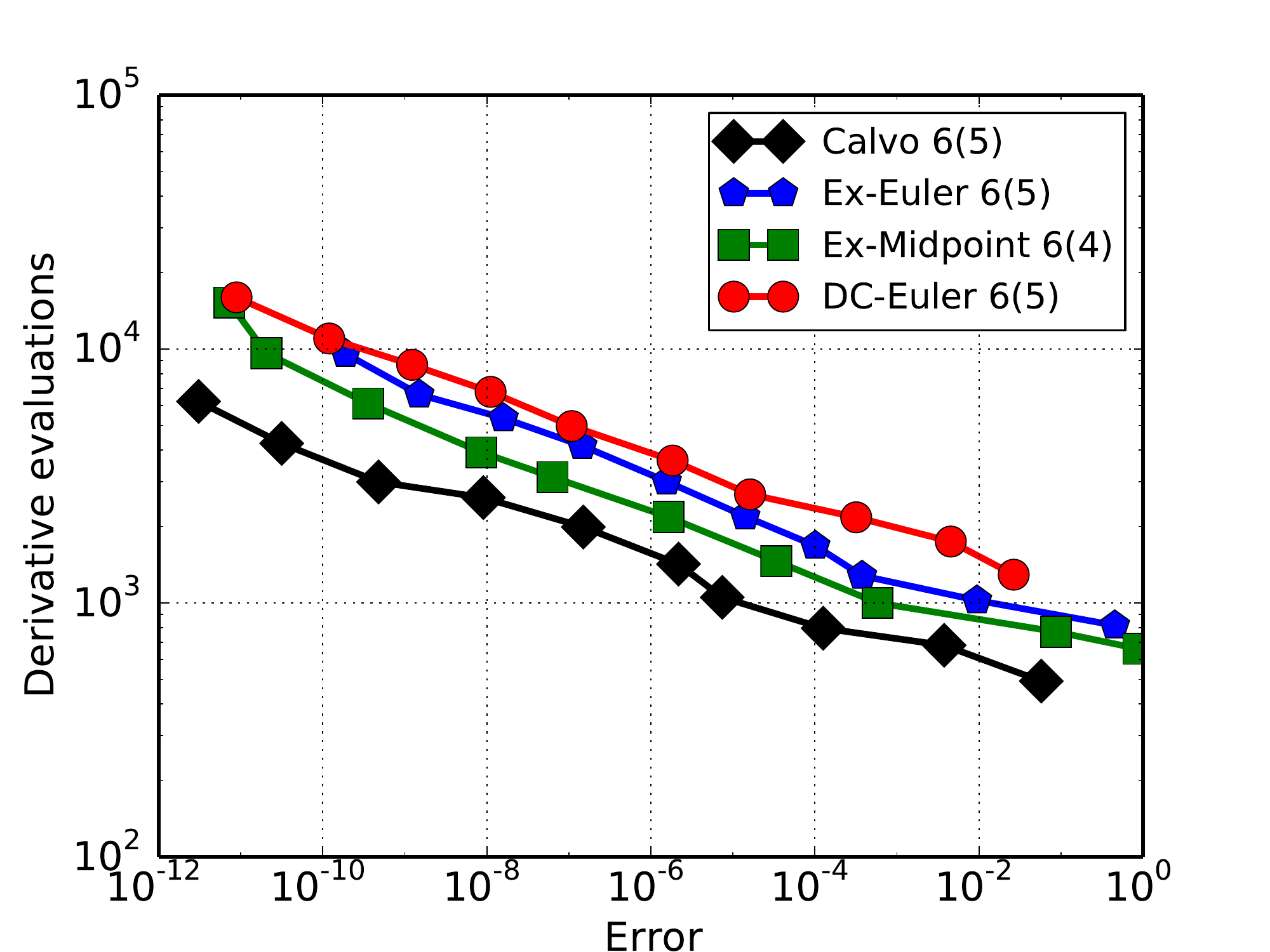} }
\subfigure[8th order]{\label{fig:d1_8}
\includegraphics[width=0.4\textwidth,height=0.32\textwidth]{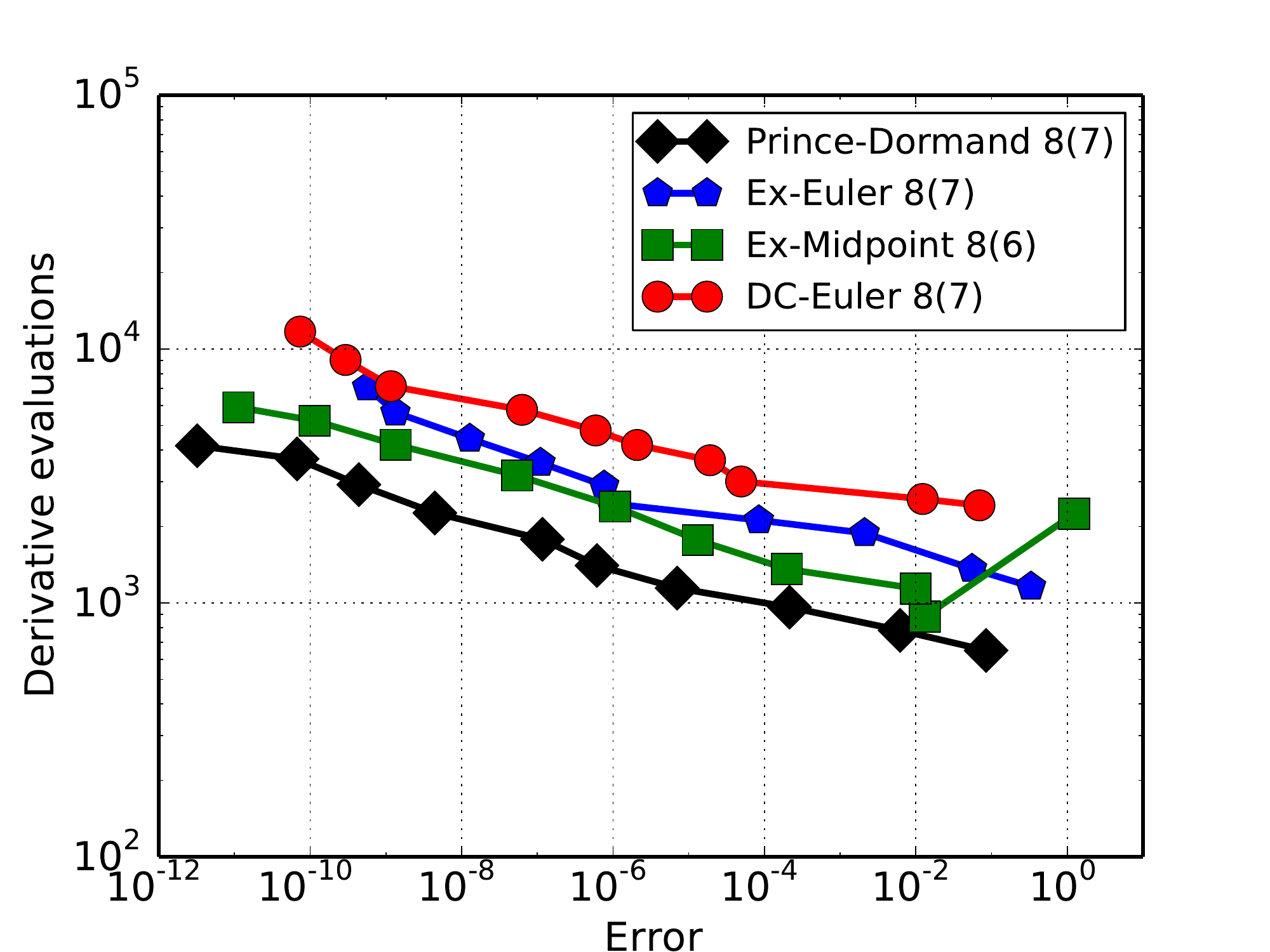} } \\
\subfigure[10th order]{\label{fig:d1_10}
\includegraphics[width=0.4\textwidth,height=0.32\textwidth]{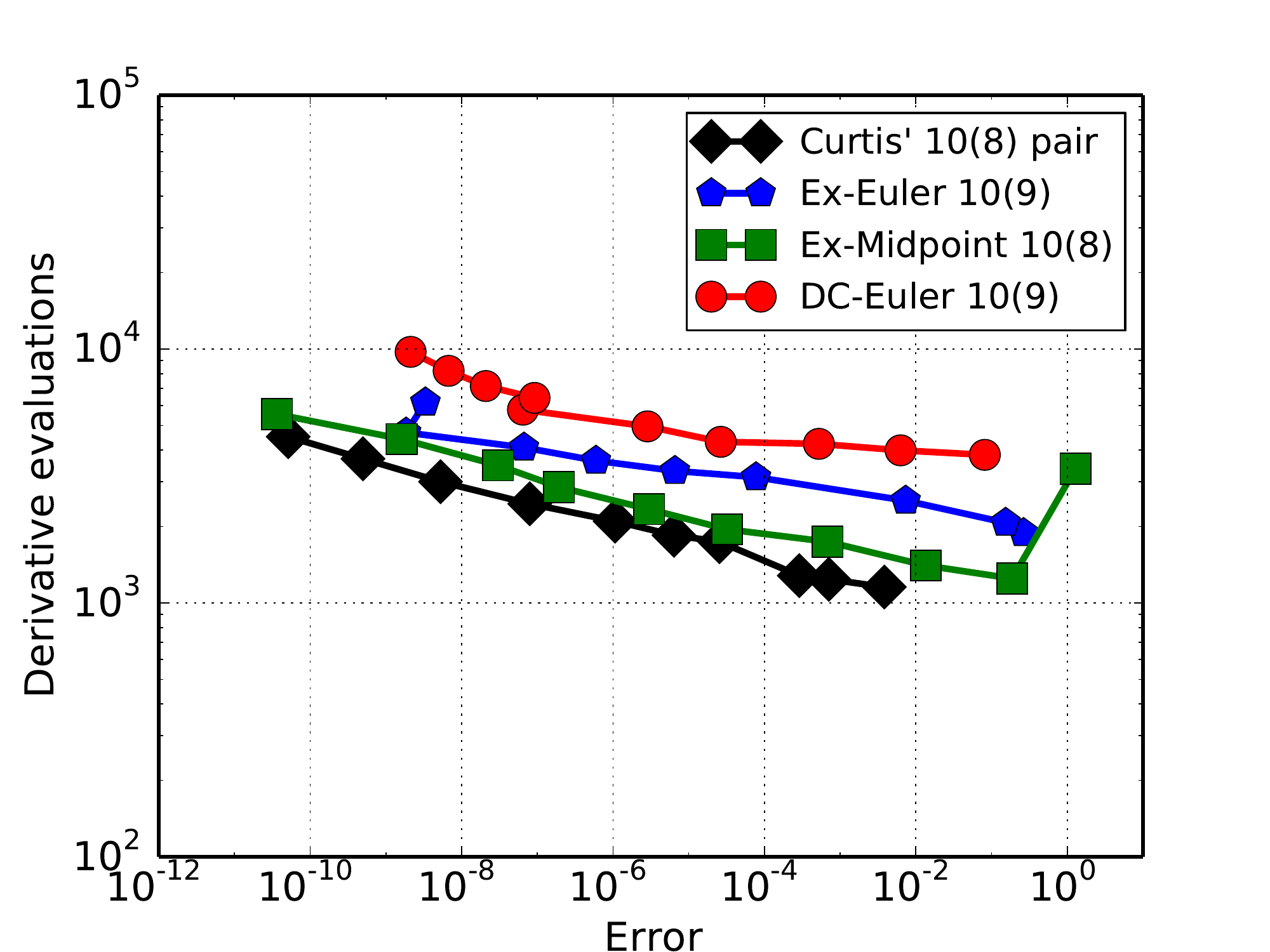} }
\subfigure[12th order]{\label{fig:d1_12}
\includegraphics[width=0.4\textwidth,height=0.32\textwidth]{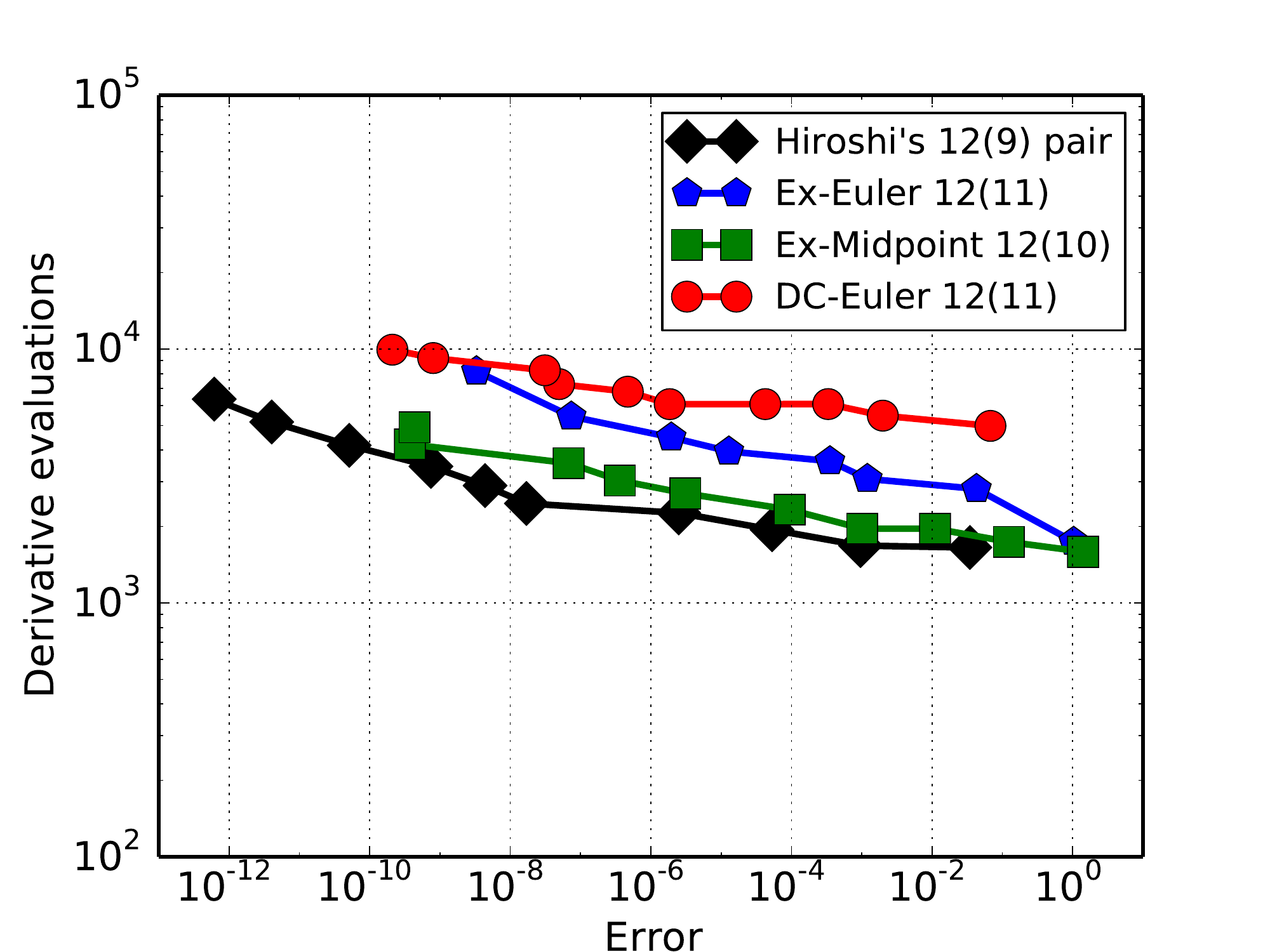} }
\end{center}
\caption{Efficiency tests on problem SB1 (Section \ref{threebody}).\label{fig:test1}}
\end{figure}

We consider the first three-body problem from~\cite{Shampine1986}:
\begin{equation}
\begin{aligned}
\text{SB1}: \qquad \qquad y_1' &= y_3,\\
y_2' &= y_4,\\
y_3' &= y_1+2y_4-\mu'\frac{y_1+\mu}{\left((y_1+\mu)^2+y_2^2\right)^{\frac{3}{2}}}-\mu\frac{y_1-\mu'}{\left((y_1-\mu')^2+y_2^2\right)^{\frac{3}{2}}}\\
y_4' &= y_2+2y_3-\mu'\frac{y_2}{\left((y_1+\mu)^2+y_2^2\right)^{\frac{3}{2}}}-\mu\frac{y_2}{\left((y_1-\mu')^2+y_2^2\right)^{\frac{3}{2}}},
\end{aligned}
\label{eq:SB1}
\end{equation}
Here $\mu'=1-\mu$, the final time is $T=6.192169331319639$, and the initial values are
\begin{align}
    y_1(0)=1.2,\;y_2(0)=0,\; y_3(0)=0,\; y_4(0)=-1.049357509830319 \text{ and } \mu=0.0121285627653123.
\end{align}
Figure~\ref{fig:test1} plots number of function evaluations (cost) against the
absolute error for this problem.  The absolute error is
\begin{equation}
\textup{Error} = |y_N-y(T)|,
\end{equation}
where $T$ is the final time and $y_N$ is the numerical solution at that time,
while $y(T)$ is a reference solution computed using a fine grid and the method
of Bogacki \& Shampine~\cite{Bogacki1996}.
The initial step size is 0.01.  In every case, the method
efficiencies follow the ordering predicted by the accuracy efficiency index,
and are consistent with previous studies.

\subsubsection{A two-population growth model\label{popgrowth}}
Next we consider problem B1 of \cite{Hull1972}, which models the growth of
two conflicting populations:
\begin{subequations}
\begin{align}
    y_1' & = 2(y_1 - y_1 y_2) & y_1(0) & = 1 \\
    y_2' & = -(y_2-y_1 y_2) & y_2(0) & 3.
\end{align}
\end{subequations}

\begin{figure}
\begin{center}
\subfigure[6th order]{\label{fig:b1_6}
\includegraphics[width=0.4\textwidth,height=0.32\textwidth]{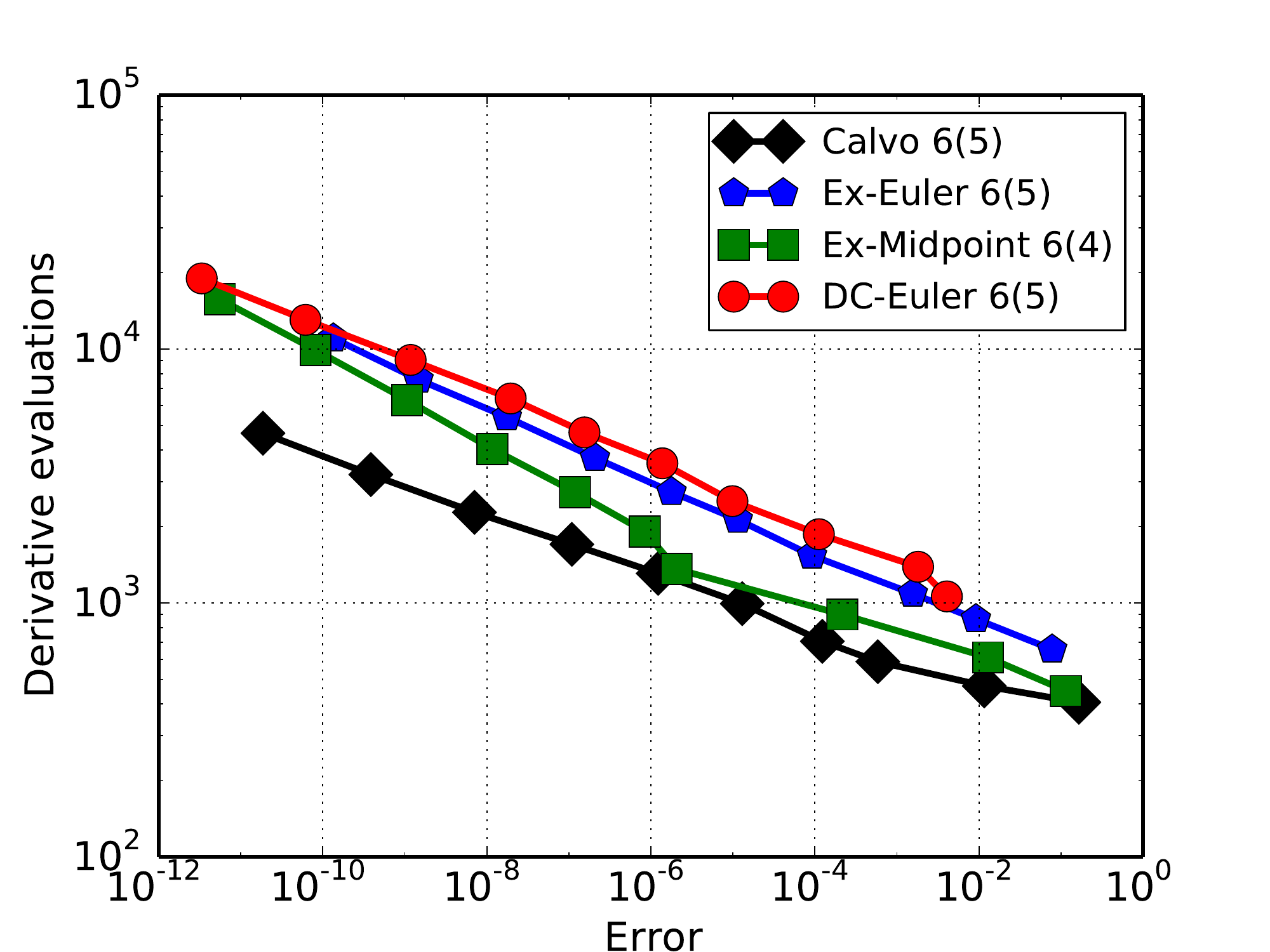} }
\subfigure[8th order]{\label{fig:b1_8}
\includegraphics[width=0.4\textwidth,height=0.32\textwidth]{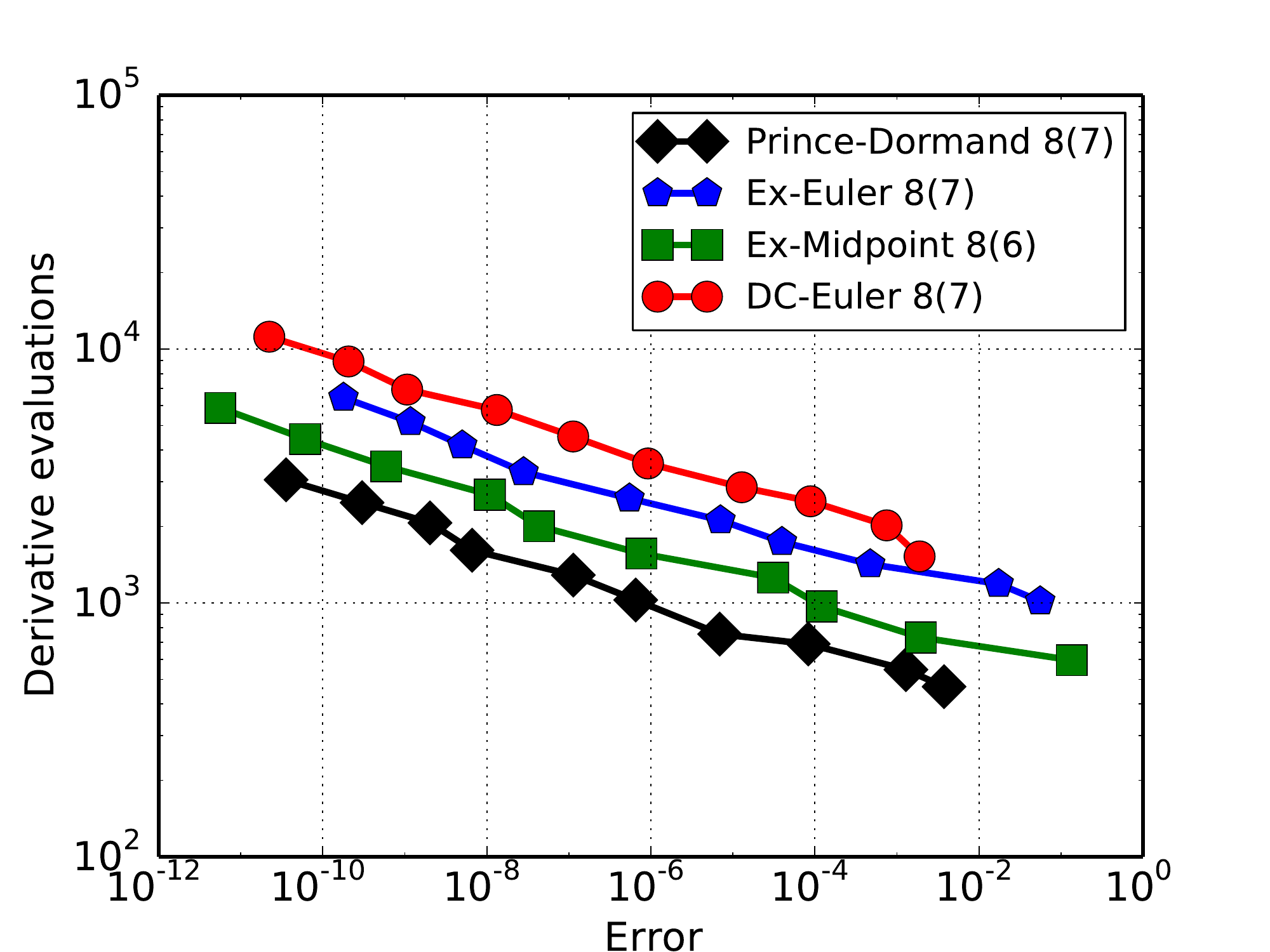} } \\
\subfigure[10th order]{\label{fig:b1_10}
\includegraphics[width=0.4\textwidth,height=0.32\textwidth]{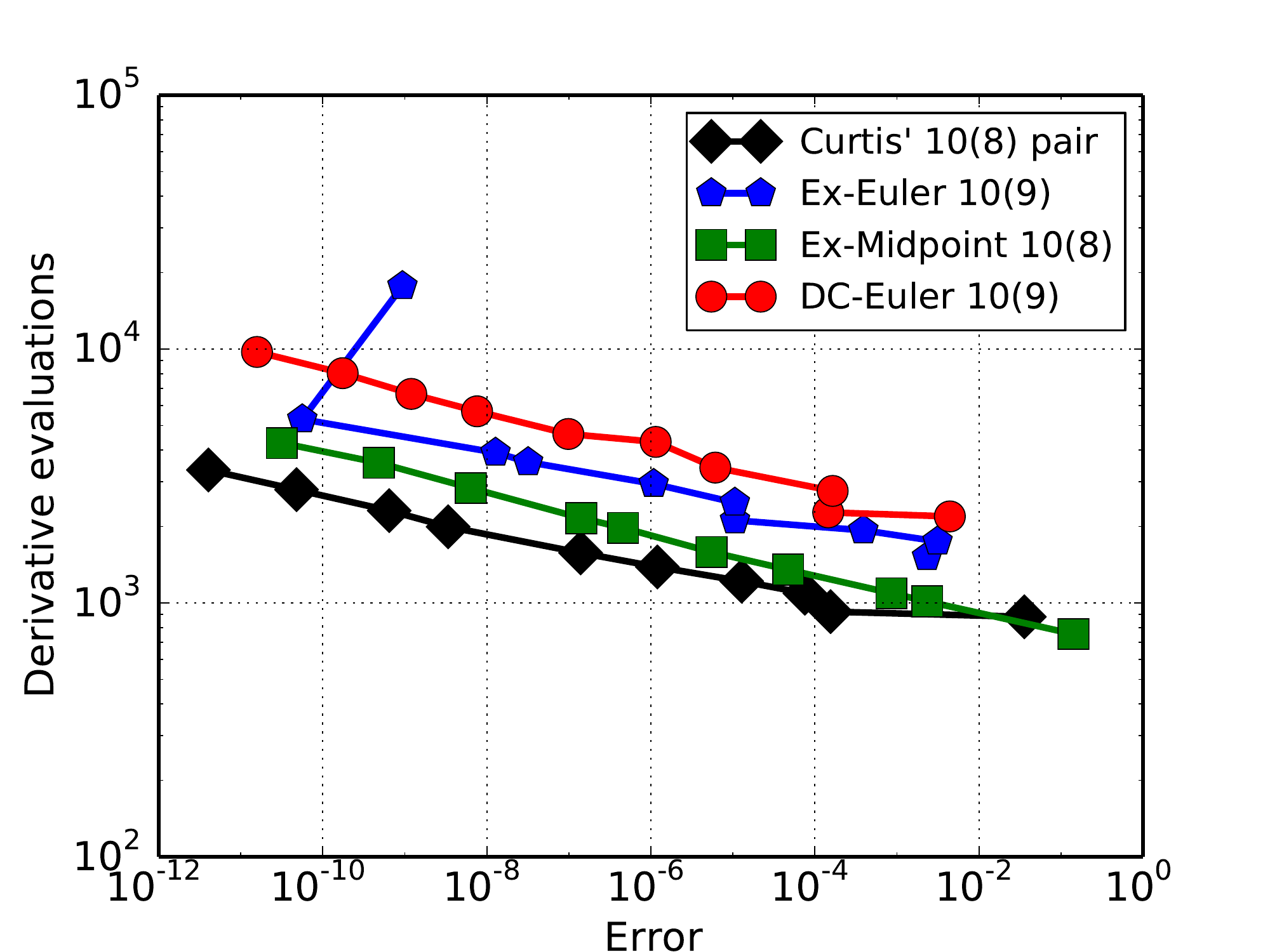} }
\subfigure[12th order]{\label{fig:b1_12}
\includegraphics[width=0.4\textwidth,height=0.32\textwidth]{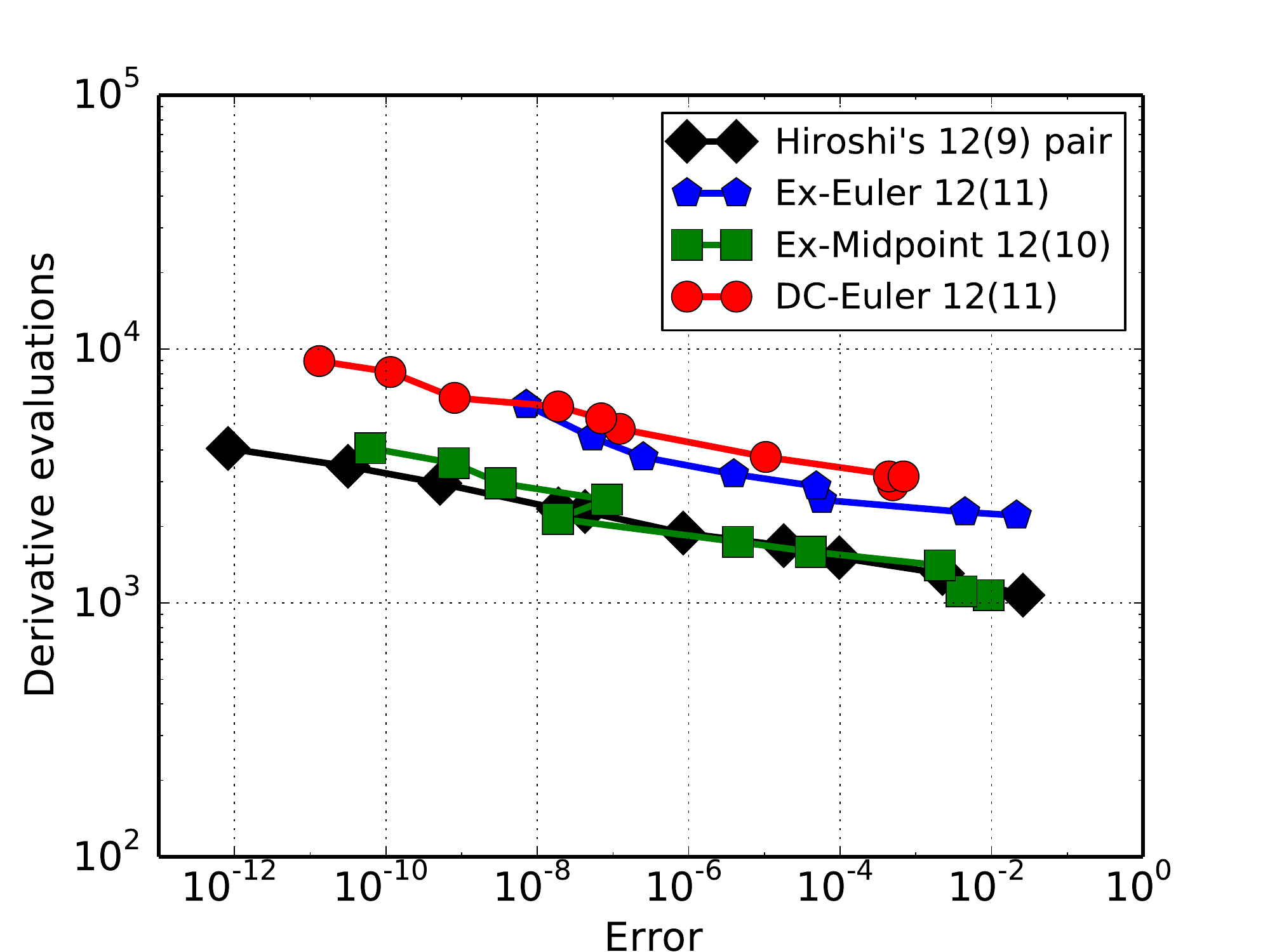} }
\end{center}
\caption{Efficiency tests on problem B1 (Section \ref{popgrowth}).\label{fig:test2}}
\end{figure}

\subsubsection{A nonlinear wave PDE\label{stegoton}}
\begin{figure}
\begin{center}
\subfigure[Solution (stress)]{\label{fig:stego}
\includegraphics[width=0.4\textwidth]{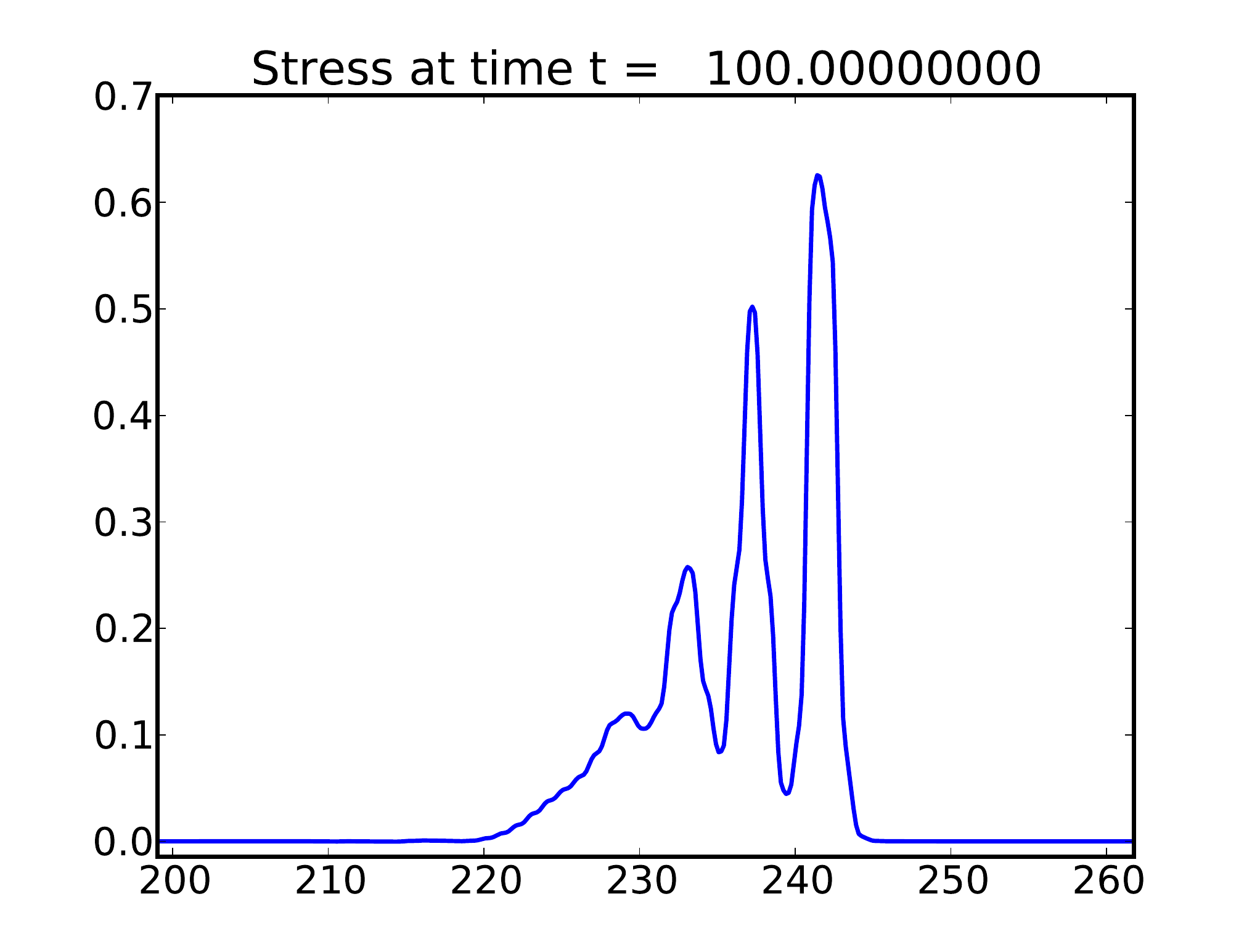} }
\subfigure[Efficiency]{\label{fig:stego_eff}
\includegraphics[width=0.4\textwidth]{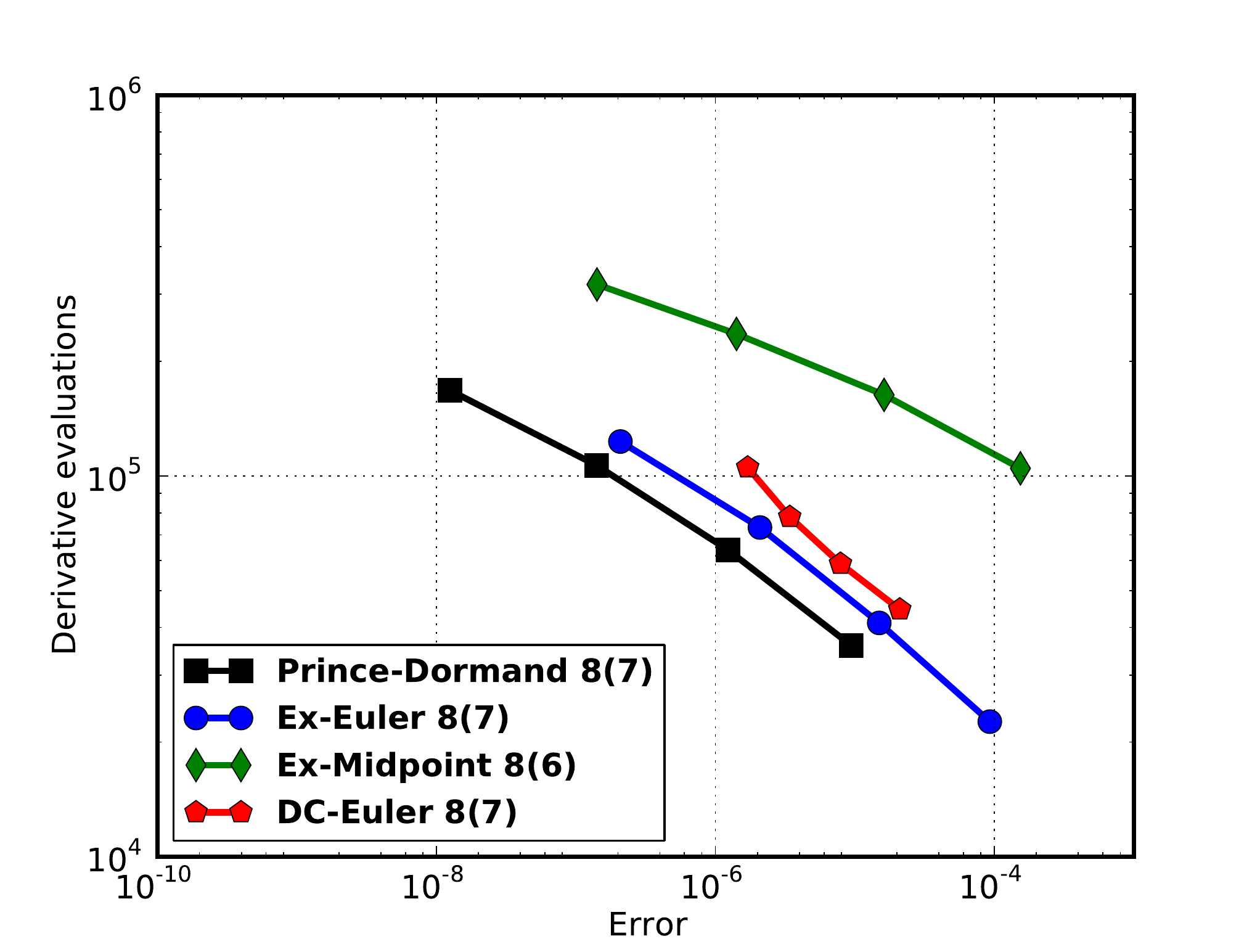} }
\end{center}
\caption{%
    Solution and efficiency for methods applied to the stegoton problem (Section \ref{stegoton}).
}%
\label{fig:stego_full}
\end{figure}

Finally, we consider the integration of a high-order PDE semi-discretization from 
\cite{ketchesonleveque_periodic}.  We solve the 1D elasticity equations
\begin{subequations}
\begin{align}
\epsilon_t(x,t)-u_x(x,t) & = 0 \\
\rho(x)u(x,t)_t - \sigma(\epsilon(x,t),x)_x & = 0.
\end{align}
\end{subequations}
with nonlinear stress-strain relation
\begin{align}
\sigma(\epsilon,x) = \exp(K(x)\epsilon) - 1,
\end{align}
and a simple periodic medium
composed of alternating homogeneous layers:
\begin{align}
\rho(x)=K(x) & = \left\{ \begin{array}{ll}
 4 & \text{if } j< x < (j+1/2)
   \mbox{ for some integer j}, \\
 1 & \mbox{otherwise.} \end{array}\right.
\end{align}
We consider the domain $0\le x \le 300$, an initial Gaussian perturbation to
the stress, and final time $T=100$.  The solution consists of two trains of 
emerging solitary waves; one of them is depicted in Figure~\ref{fig:stego}.
The semi-discretization is based on the WENO wave-propagation method implemented
in SharpClaw~\cite{Ketcheson2011}.

Efficiency results for 8th order methods are shown in Figure \ref{fig:stego_eff}.
The spatial grid is held fixed across all runs, and the time step is adjusted
automatically to satisfy the imposed tolerance.
The error is computed with respect to a solution computed with tolerance $10^{-13}$
using the 5(4) pair of Bogacki and Shampine.
For the most part, these are consistent with the results from the smaller problems
above.  However, the midpoint extrapolation method performs quite poorly on this
problem.  The reason is not clear, but this underscores the fact that performance
on particular problems can be very different from the ``average'' performance of a method.

\subsection{Failure of integrators\label{sec:failure}}
Some failure of the integrators was observed in testing.  These failures fall 
into two categories.  First, at very tight tolerances, the high order Euler
extrapolation methods were sometimes unable to finish because the time step size was driven
to zero.  This is a known issue related to internal stability; for a full
explanation see \cite{klp_internal}.

Second, the deferred correction methods sometimes gave
global errors much larger than those obtained with the other methods.  This
indicates a failure of the error estimator.  Upon further investigation, we
found that the natural embedded error estimator method of order $p-1$ satisfies
nearly all (typically all but one) of the order conditions for order $p$.
Hence these estimators may be said to be {\em defective}, and it would
be advisable to employ a more robust approach like that discussed in \cite{dutt2000}.
Since our focus is purely on Runge--Kutta pairs, we do not pursue this issue
further here.

\subsection{Ideal parallel performance\label{sec:odex}}
\begin{figure}
\begin{center}
\subfigure[6th order methods]{
\includegraphics[width=0.4\textwidth]{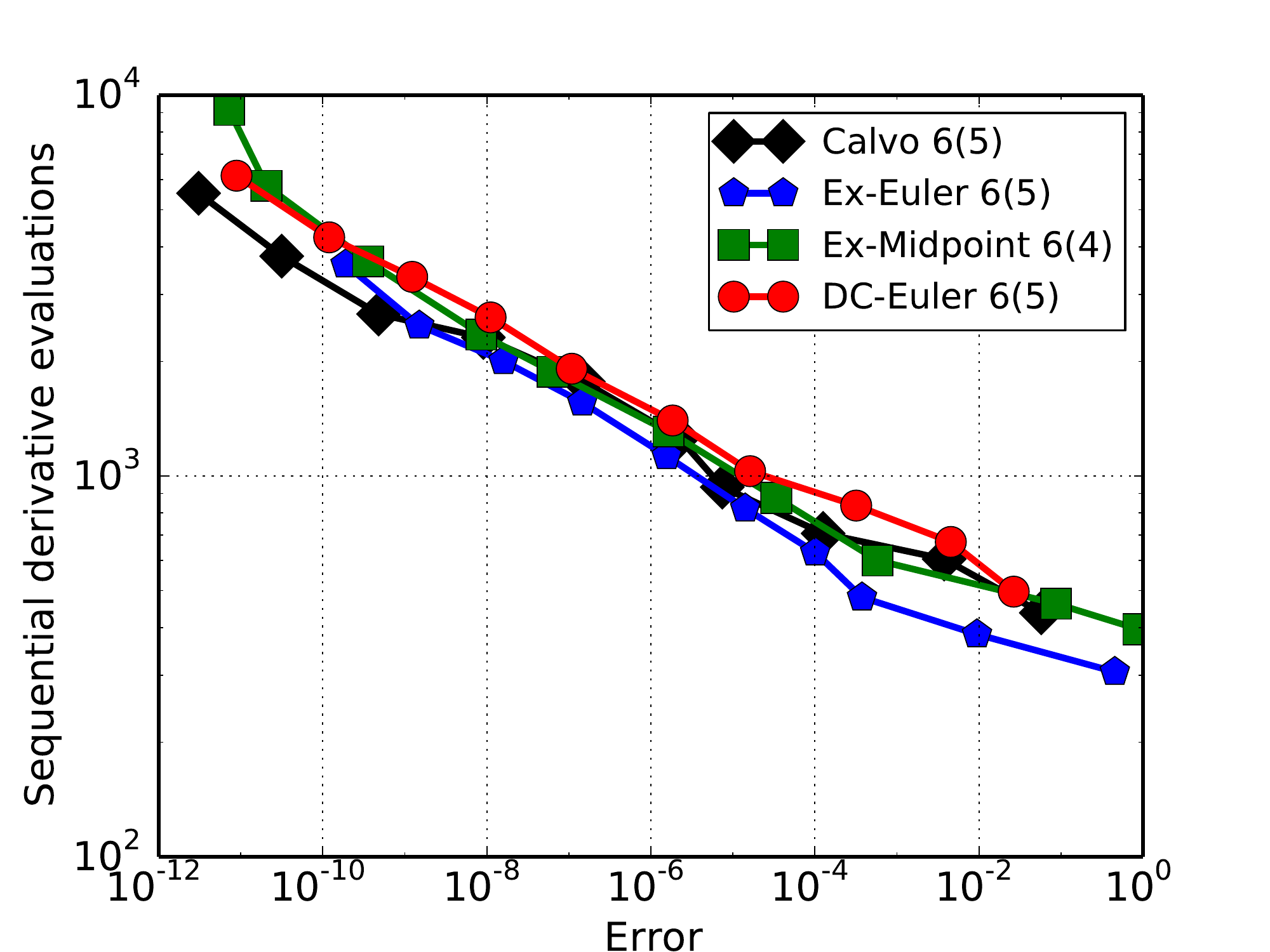}}
\subfigure[8th order methods]{
\includegraphics[width=0.4\textwidth]{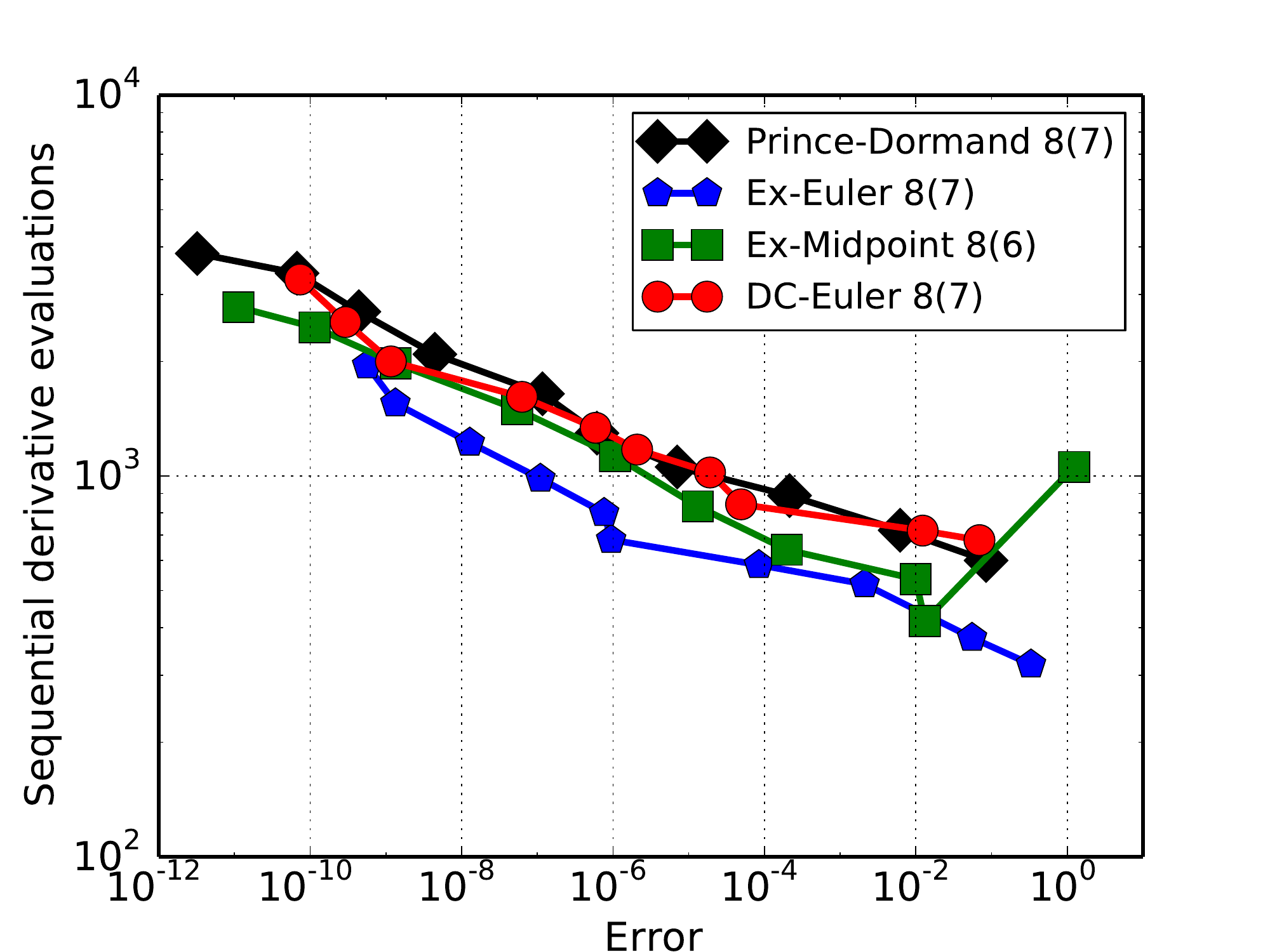}}
\subfigure[10th order methods]{
\includegraphics[width=0.4\textwidth]{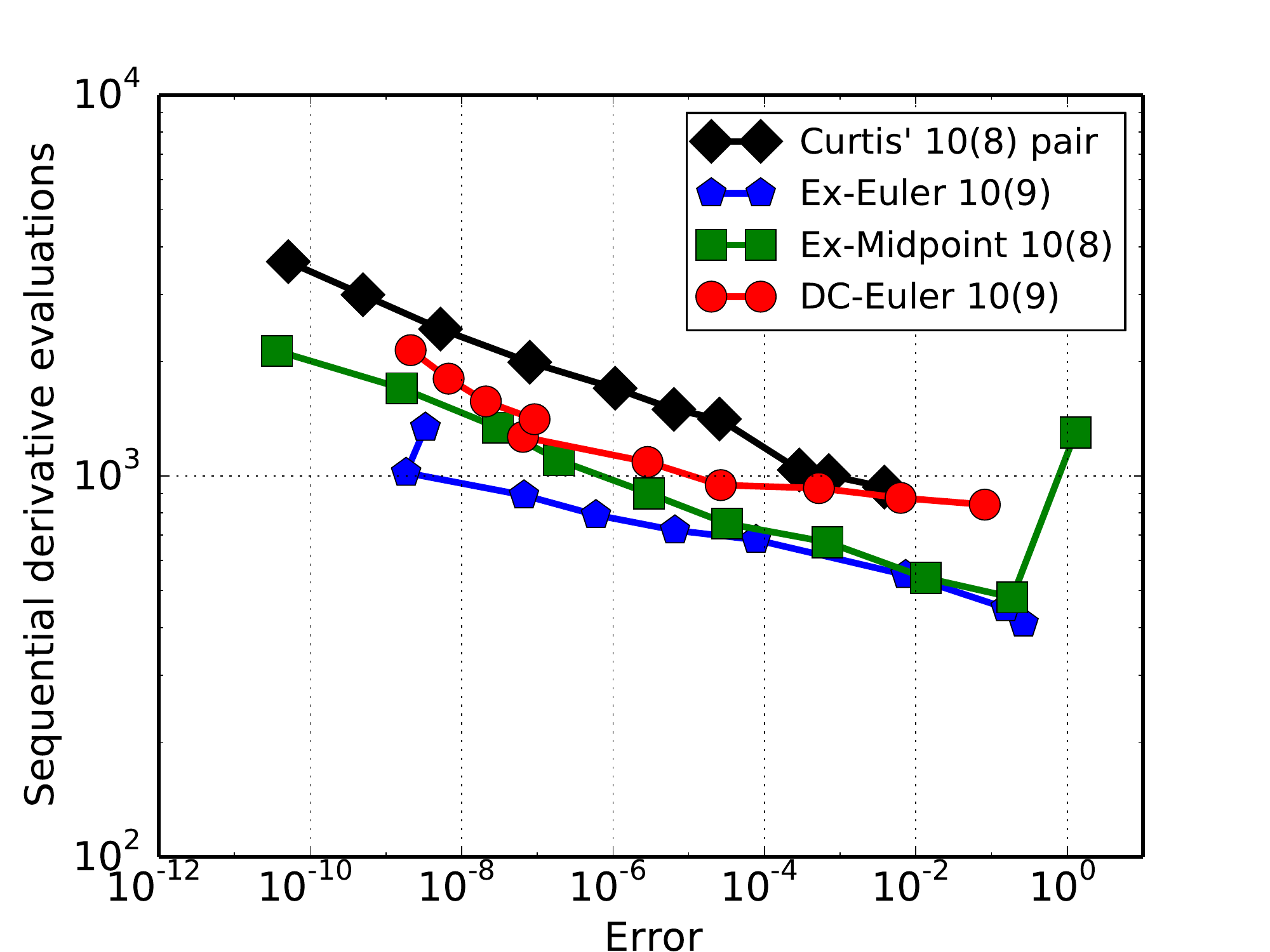}}
\subfigure[12th order methods]{
\includegraphics[width=0.4\textwidth]{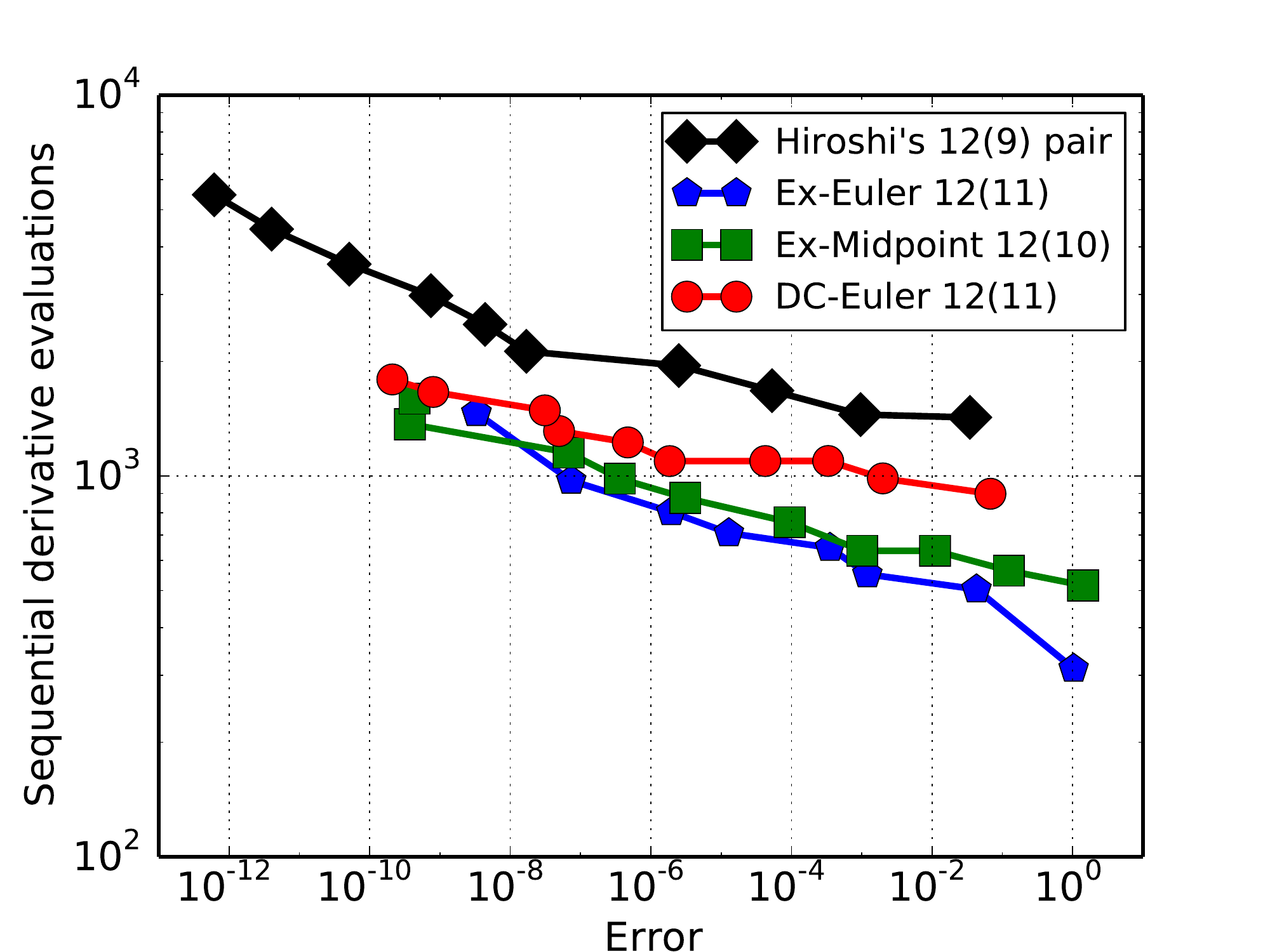}}
\caption {Efficiency for the problem SB1 (Section \ref{threebody}) based on sequential derivative evaluations.}
\label{fig:sb1_par}
\end{center}
\end{figure}

\begin{figure}
\begin{center}
\subfigure[6th order methods]{
\includegraphics[width=0.4\textwidth]{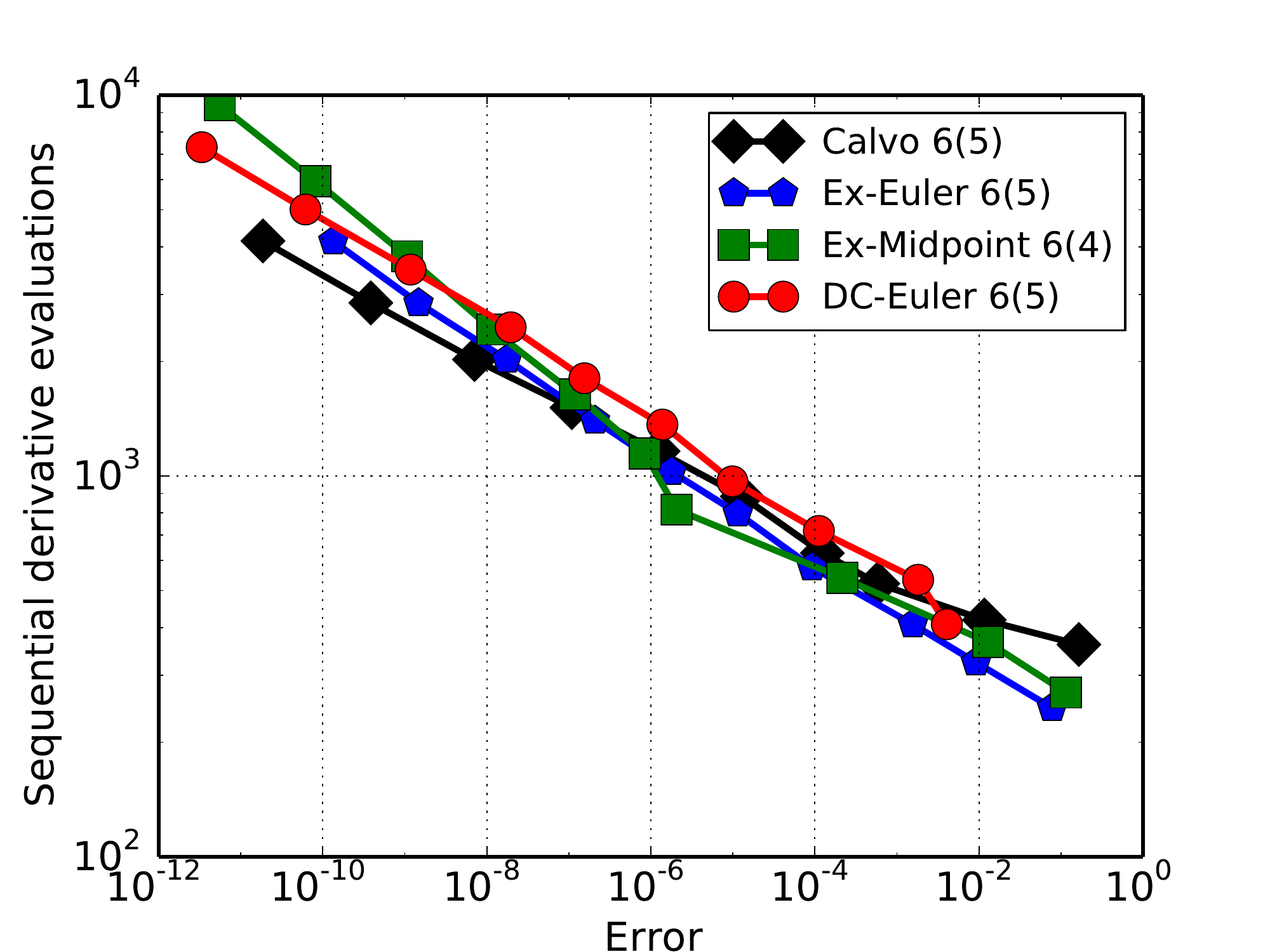}}
\subfigure[8th order methods]{
\includegraphics[width=0.4\textwidth]{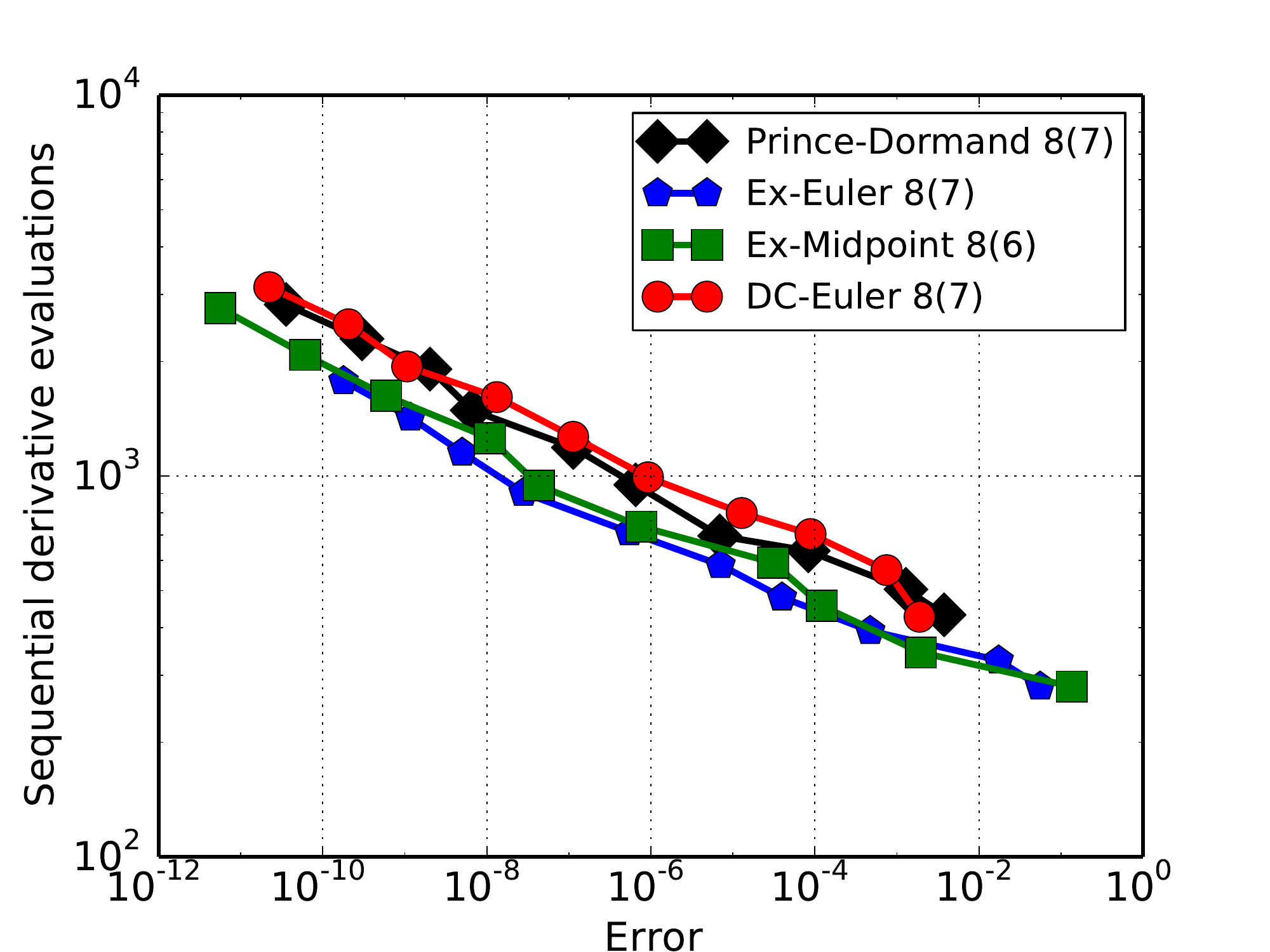}}
\subfigure[10th order methods]{
\includegraphics[width=0.4\textwidth]{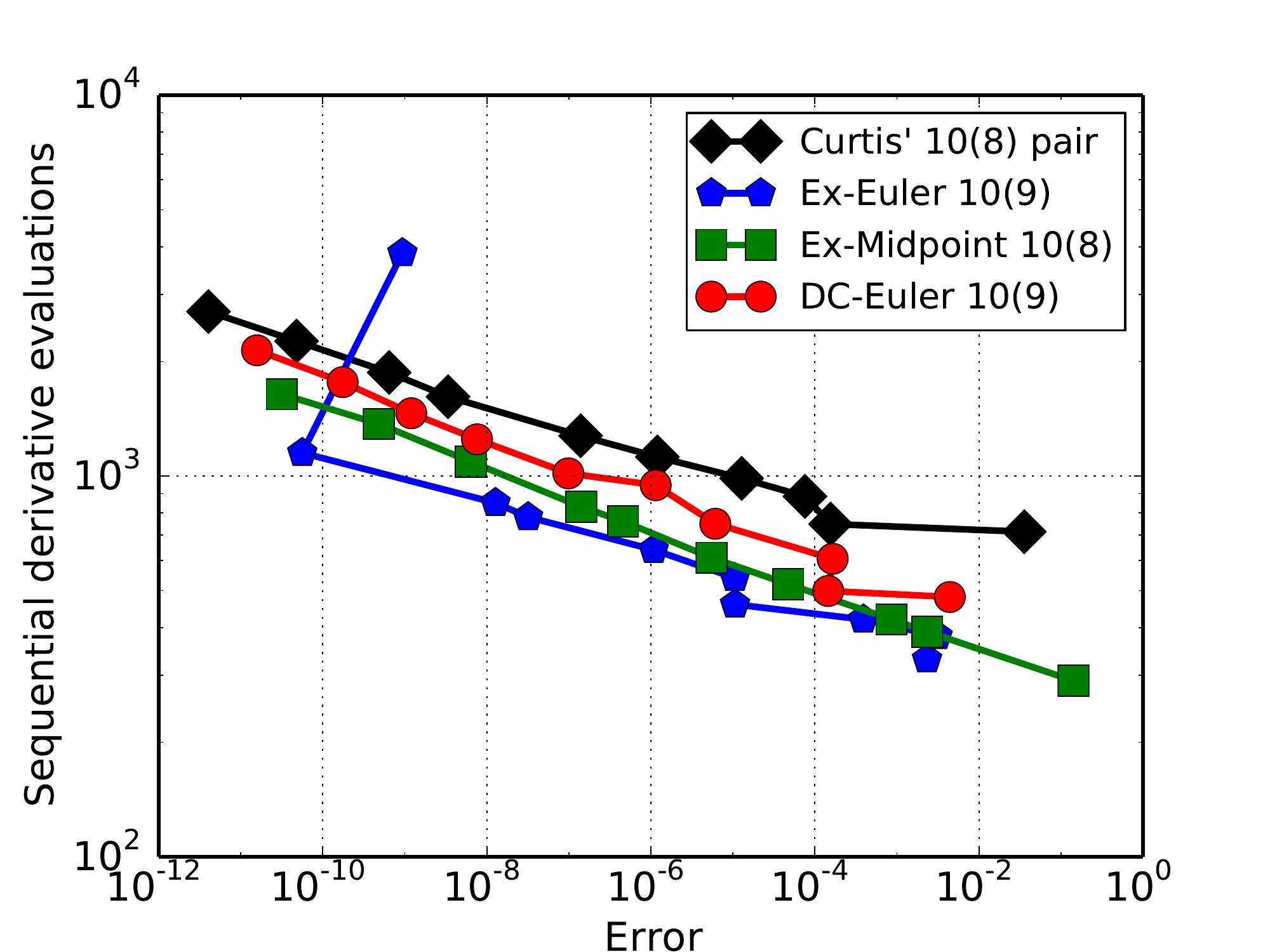}}
\subfigure[12th order methods]{
\includegraphics[width=0.4\textwidth]{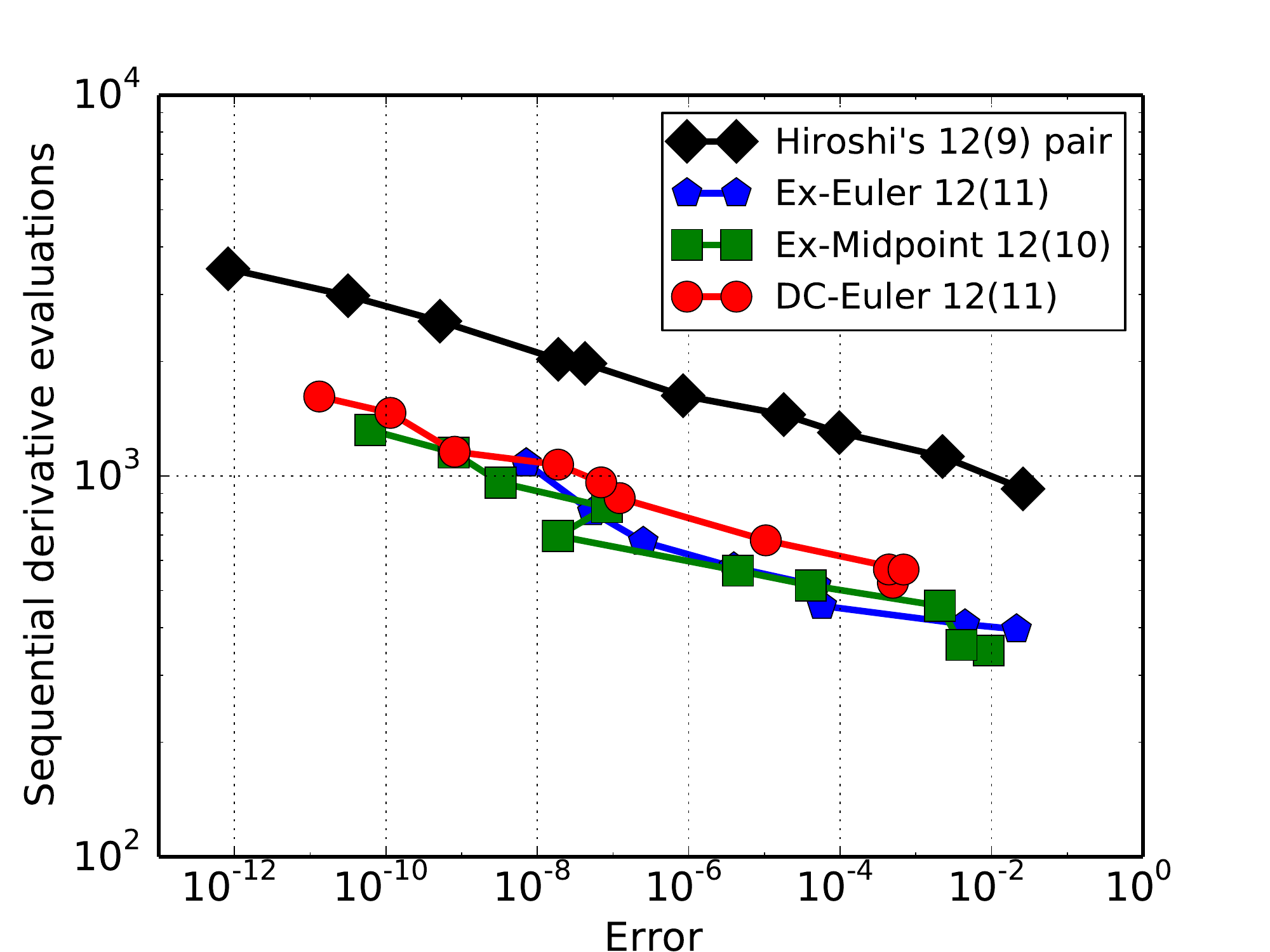}}
\caption {Efficiency for the problem B1 (Section \ref{popgrowth}) based on sequential derivative evaluations.}
\label{fig:b1_par}
\end{center}
\end{figure}

\begin{figure}
\begin{center}
\includegraphics[width=0.4\textwidth]{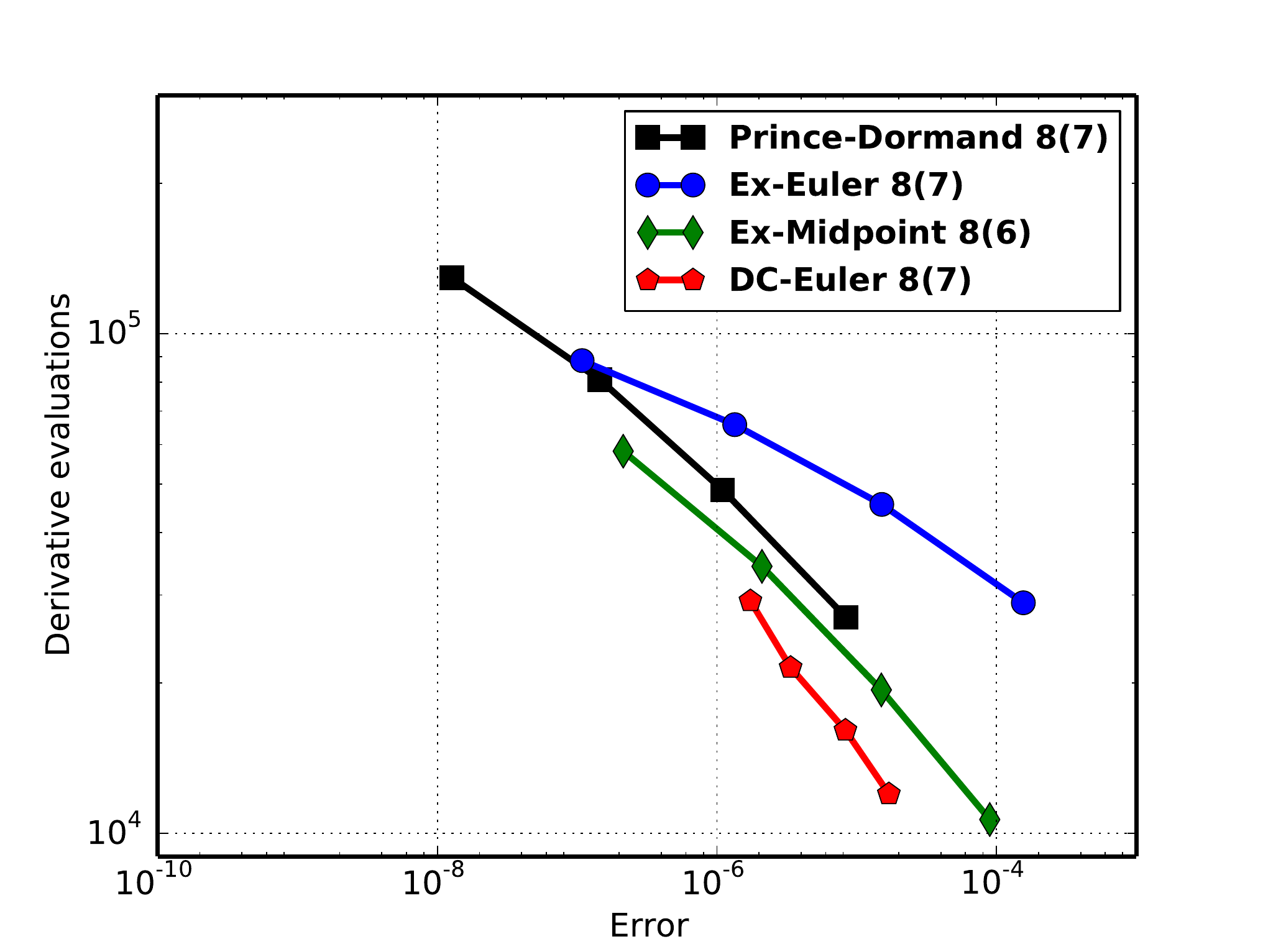}
\end{center}
\caption{%
    Parallel efficiency for the stegoton problem (Section \ref{stegoton}).
\label{fig:stego_eff_par}}

\end{figure}

Figures~\ref{fig:sb1_par}-\ref{fig:stego_eff_par} show efficiency for the same three test problems
but now based on the number of sequential evaluations.  That is, the vertical axis is $N\spar$, where $N$
denotes the number of steps taken.  The same measure of efficiency was used in
\cite{van1990parallel}.  We see that the parallelizable methods -- especially extrapolation --
outperform traditional methods, especially at higher orders.  Similar results
were obtained for parallel iterated RK methods in \cite{van1990parallel}.
Remarkably, the deferred correction method performs the best by this measure for the stegoton problem.

This measure of efficiency may be viewed with some skepticism since it neglects
the cost of communication.  This concern is addressed with a true parallel implementation
in the next section.

\section{A shared-memory implementation of extrapolation\label{sec:shared}}
Development and testing of a tuned
parallel extrapolation or deferred correction code is beyond the scope of this paper,
but in this section we run a simple example to demonstrate that it is possible in practice 
to achieve speedups like those listed in Table~\ref{tbl:par}, and to outperform
even the best highly-tuned traditional RK methods, at least on problems with an
expensive right-hand-side.
We focus on speedup with an eye to providing efficient black-box
parallel ODE integrators for multicore machines, noting that the number of available
cores is often more than can be advantageously used by the methods considered.

Previous studies have implemented explicit extrapolation methods in parallel
and achieved parallel efficiencies of up to about 80\%
\cite{Ito_Fukushima_1997,Kappeller_Kiehl_Perzl_Lenke_1996,lustman1992solution}.
As those studies were conducted about twenty years ago, it is not clear
that their conclusions are relevant to current hardware.

In order to test the achievable parallel speedup, we took the code ODEX \cite{hairer1993},
downloaded from \url{http://www.unige.ch/~hairer/software.html}, and
modified it as follows:
\begin{itemize}
    \item Fixed the order of accuracy (disabling adaptive order selection)
    \item Inserted an OMP PARALLEL pragma around the extrapolation loop
    \item Removed the smoothing step
\end{itemize}
We refer to the modified code as ODEX-P.

\begin{figure}
\begin{center}
\subfigure[Dynamic scheduling\label{fig:naive}]{
\includegraphics[width=0.4\textwidth]{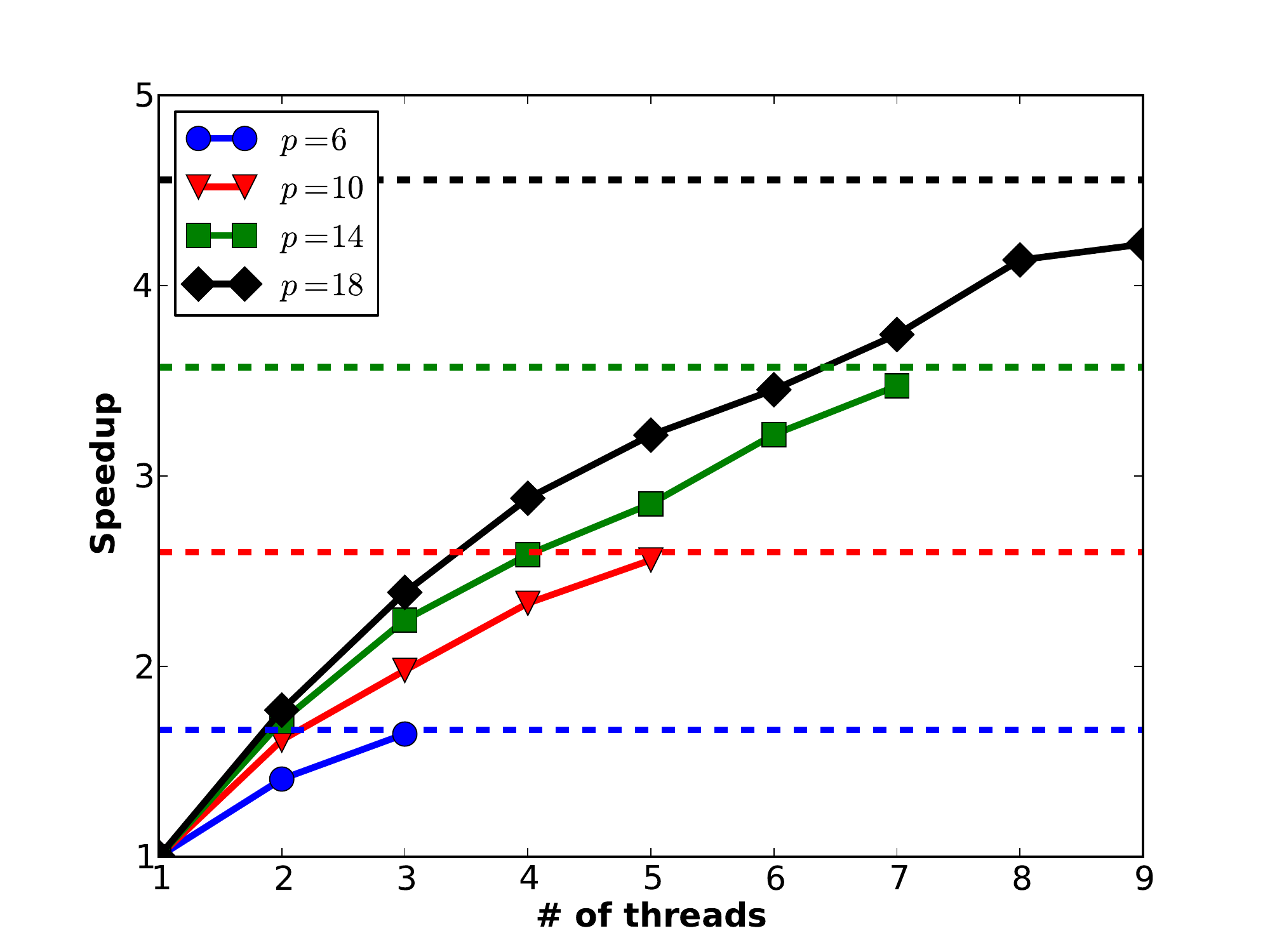}}
\subfigure[Manual load-balancing\label{fig:load_balanced}]{
\includegraphics[width=0.4\textwidth]{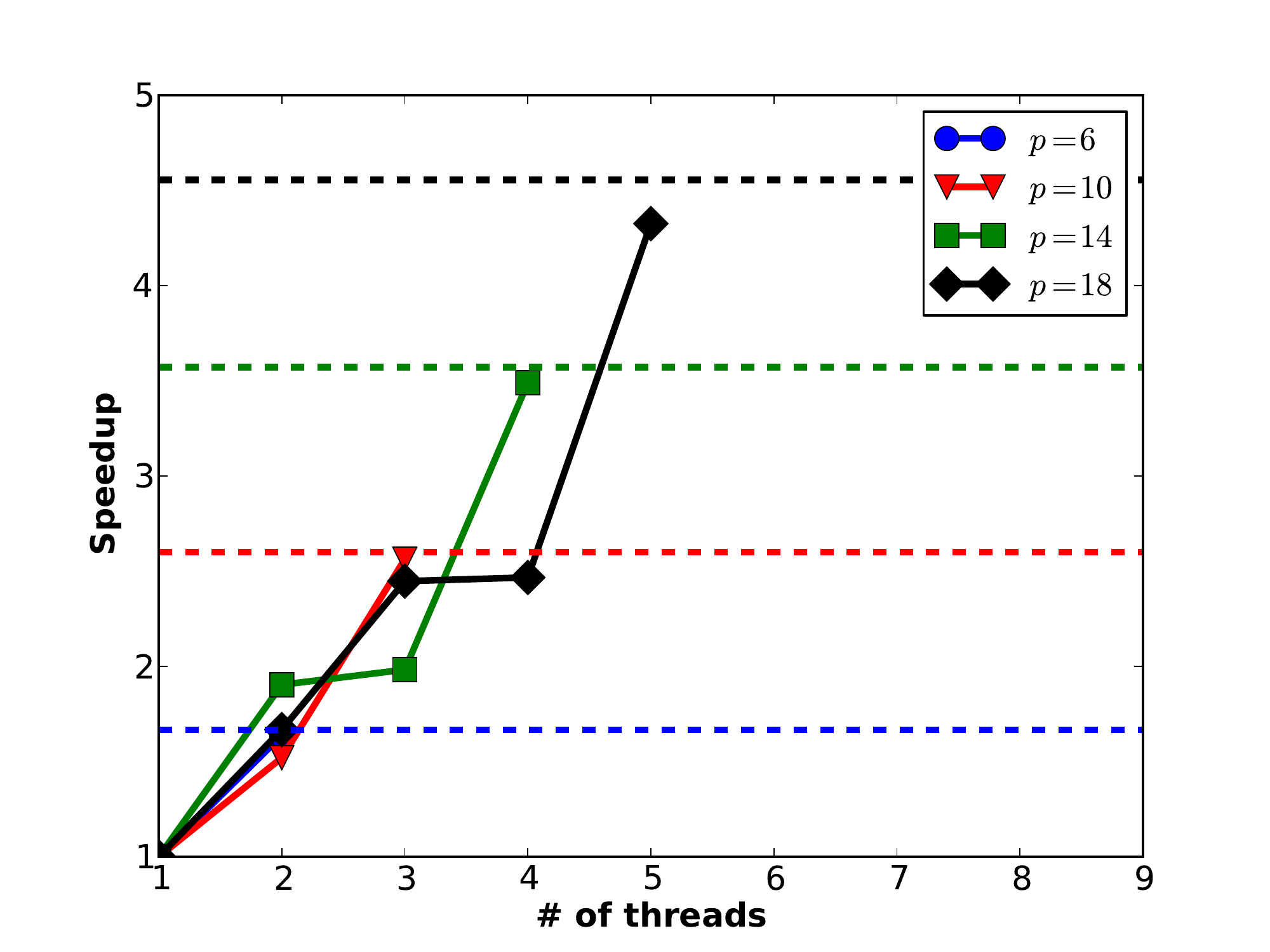}}
\caption{Measured speedup of the midpoint extrapolation code ODEX-P on a 400-body
gravitation problem by insertion of a single OMP parallel pragma in the code.
The ratio of runtime with multiple threads to runtime using a single thread is
plotted.
The dotted lines show the theoretical
maximum speedup $S = (p^2+4)/(4p)$ based on our earlier analysis.
\label{fig:odex_speedup}}
\end{center}
\end{figure}

Figure \ref{fig:naive} shows the achieved speedup based on
dynamic scheduling for $p=6,10,14,18$,
applying the code to an $N$-body gravitational problem with 400 bodies.
Results for other orders are similar.  The dotted lines show the theoretical
maximum speedup $S = (p^2+4)/(4p)$ based on our earlier analysis.
The tests were run on a workstation with two 2.66 Ghz quad-core Intel Xeon processors,
and the code was compiled using gfortran.
Using $p/2$ threads, the measured speedup is very close to the theoretical
maximum.  However, the speedup is significantly below the theoretical value
when only $P$ threads are used.  We interpret this to mean that the dynamic
scheduler is not able to optimally allocate the work among threads unless
there are enough threads to give just one loop iteration to each.

Figure \ref{fig:load_balanced} and Table \ref{tbl:odex}
show the result of a more intelligent
parallel implementation, using static scheduling with the code modified so that
both $T_{k1}$ and $T_{r-k,1}$ are computed in a single loop iteration.  
This load balancing scheme is optimal when using on the optimal number of
threads $P$, and the results agree almost perfectly with theory.

\begin{table}
\centering
\begin{tabular}{c|c|c|c|c|c|c|c}
& & \multicolumn{2}{c}{Runtime} & \multicolumn{2}{|c}{Max. speedup} & \multicolumn{2}{|c}{Parallel Efficiency} \\
Order ($p$) & $P$ & 1 thread & $P$ threads & Theory ($S$) & Observed & Theory ($E$) & Observed \\ \hline
6           & 2 &   13.140 & 7.977       & 1.67 &  1.65        & 0.83  & 0.82 \\ \hline
10          & 3 &   17.370 & 6.770       & 2.60 &  2.57        & 0.87  & 0.86 \\ \hline
14          & 4 &   19.508 & 5.573       & 3.57 &  3.50        & 0.89  & 0.88 \\ \hline
18          & 5 &   25.876 & 5.827       & 4.56 &  4.44        & 0.91  & 0.89 \\ \hline
\end{tabular}
\caption{Runtime, speedup and efficiency of manually load-balanced runs of the 
modified ODEX-P code with $P$ threads.  The observed speedup (and efficiency) are close
to the theoretically optimal values ($S$ and $E$).\label{tbl:odex}}
\end{table}

\subsection{Comparison with DOP853\label{sec:runtimes}}
We now compare actual runtimes of our experimental ODEX-P with
the DOP853 code available from \url{http://www.unige.ch/~hairer/prog/nonstiff/dop853.f}.
These two codes have been compared in \cite[Section II.10]{hairer1993}, but using the original
ODEX code (with order-adaptivity and without parallelism).  In that reference, DOP853 was
shown to be superior to ODEX at all but the most strict tolerances.

\begin{figure}
\begin{center}
\includegraphics[width=0.5\textwidth]{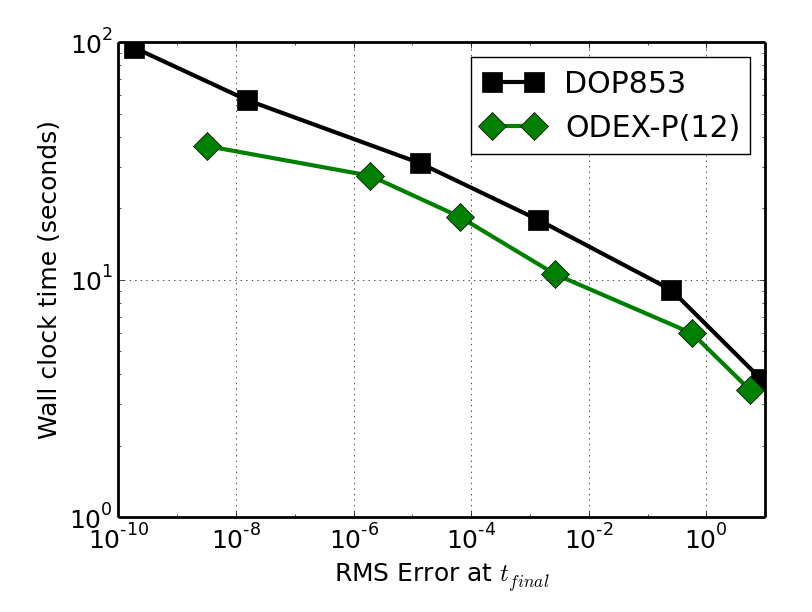}
\caption{Runtime versus achieved relative error of the midpoint extrapolation
code ODEX-P on a 400-body gravitation problem.
The tests were run on a workstation with two 2.66 Ghz quad-core Intel Xeon processors
using four threads.
\label{fig:realtimes}}
\end{center}
\end{figure}

Table~\ref{tbl:realtimes} shows runtimes versus prescribed tolerance
for a 400-body problem for the two codes, using fixed order 12 (with 4 threads)
in the ODEX-P code.  Figure~\ref{fig:realtimes} shows the achieved relative
root-mean-square global error versus runtime.
Perhaps surprisingly, the parallel extrapolation
code is no worse even at loose tolerances.  At moderate to strict tolerances, it 
substantially outperforms the RK code.

\begin{table}
\centering
\begin{tabular}{l|c|c|c|c|c}
& \multicolumn{5}{c}{Tolerance} \\
Code       & 1.e-3 & 1.e-5 & 1.e-7 & 1.e-9 & 1.e-11 \\ \hline \hline
DOP853     & 3.81  &  9.02 & 17.80 & 30.93 & 56.75  \\ \hline
ODEX-P(12) & 3.31  &  5.87 & 10.32 & 18.07 & 26.76
\end{tabular}
\caption{Runtimes (in seconds) for Dormand-Prince and modified 12th-order ODEX-P code.
The tests were run on a workstation with two 2.66 Ghz quad-core Intel Xeon processors
using four threads.
\label{tbl:realtimes}}
\end{table}

\section{Discussion}
This study is intended to provide a broadly useful characterization of the properties
of explicit extrapolation and spectral deferred correction methods.
Of course, no study like this can be exhaustive.  Our approach handicaps extrapolation
and deferred correction methods by fixing the order throughout each computation;
practical implementations are order-adaptive and should achieve somewhat better efficiency.
We have investigated only the most generic versions of each class of methods;
other approaches (e.g., using higher order building blocks or exploiting
concurrency in different ways) may give significantly different results. 
Such approaches could be evaluated using the same kind of analysis employed here.
Finally, our parallel computational model is valid only when evaluation of $f$
is relatively expensive -- but that is when efficiency and concurrency
are of most interest.

The most interesting new conclusions from the present study is that
parallel extrapolation methods of very high order outperform
sophisticated implementations of the best available RK methods for problems
with an expensive right hand side.
This is true even for a relatively naive non-order-adaptive code.
We have shown that near-optimal speedup can be achieved in practice with simple
modification of an existing code.  The resulting algorithm is faster
(at least for some problems) than the highly-regarded DOP853 code.

Our serial results are in line with those of previous studies.  New here
is the evidence that spectral deferred correction -- like extrapolation --
seems inferior to well-designed RK methods (in serial).  However, we have tested
only one of the many possible variants of these methods.

High order Euler extrapolation methods suffer from dramatic amplification of roundoff errors.
This leads to the loss of several digits of accuracy (and failure of the
automatic error control) for very high order methods, and is observed in practice
on most problems.  Fortunately, midpoint extrapolation does not exhibit this
amplification.



The theoretical and preliminary experimental results we have presented suggest
that a carefully-designed parallel code based on midpoint extrapolation could be very efficient.
Such a practical implementation is the subject of current efforts.

{\bf Acknowledgment.}  We thank one of the referees, who pointed out a discrepancy that
revealed an important bug in our implementation of some spectral deferred correction methods when
$\theta\ne0$.

\bibliographystyle{plain}
\bibliography{RKextrapolation}

\begin{thebibliography}{10}

\bibitem{Bogacki1996}
P~Bogacki and Lawrence~F Shampine.
\newblock An efficient {R}unge-{K}utta (4, 5) pair.
\newblock {\em Computers \& Mathematics with Applications}, 32(6):15--28, 1996.

\bibitem{burrage1995parallel}
Kevin Burrage.
\newblock {\em Parallel and sequential methods for ordinary differential
  equations}.
\newblock Clarendon Press, 1995.

\bibitem{butcher2003}
J.~C. Butcher.
\newblock {\em {Numerical Methods for Ordinary Differential Equations}}.
\newblock Wiley, 2nd edition, 2008.

\bibitem{calvo_new_1990}
M.~Calvo, {J.I.} Montijano, and L.~Randez.
\newblock A new embedded pair of {R}unge-{K}utta formulas of orders 5 and 6.
\newblock {\em Computers \& Mathematics with Applications}, 20(1):15--24, 1990.

\bibitem{Christlieb_Macdonald_Ong_2010}
Andrew~J. Christlieb, Colin~B. Macdonald, and Benjamin~W. Ong.
\newblock Parallel high-order integrators.
\newblock {\em SIAM Journal on Scientific Computing}, 32(2):818–835, Jan
  2010.

\bibitem{curtis_high-order_1975}
A.~R. Curtis.
\newblock High-order explicit {R}unge-{K}utta formulae, their uses, and
  limitations.
\newblock {\em {IMA} Journal of Applied Mathematics}, 16(1):35–52, 1975.

\bibitem{Daniel_Pereyra_Schumaker_1968}
James~W. Daniel, Victor Pereyra, and Larry~L. Schumaker.
\newblock Iterated deferred corrections for initial value problems.
\newblock {\em Acta Cientifica Venezolana}, 19:128–135, 1968.

\bibitem{deuflhard_order_1983}
P.~Deuflhard.
\newblock Order and stepsize control in extrapolation methods.
\newblock {\em Numerische Mathematik}, 41(3):399–422, 1983.

\bibitem{Deuflhard1985}
P.~Deuflhard.
\newblock Recent progress in extrapolation methods for ordinary differential
  equations.
\newblock {\em SIAM Review}, 27(4):505--535, December 1985.

\bibitem{dutt2000}
A.~Dutt, L.~Greengard, and Vladimir Rokhlin.
\newblock {Spectral deferred correction methods for ordinary differential
  equations}.
\newblock {\em BIT Numerical Mathematics}, 40:241--266, 2000.

\bibitem{emmett2012toward}
Matthew Emmett and Michael Minion.
\newblock Toward an efficient parallel in time method for partial differential
  equations.
\newblock {\em Comm. App. Math. and Comp. Sci}, 7(1):105--132, 2012.

\bibitem{feagin2012high}
T~Feagin.
\newblock High-order explicit {R}unge-{K}utta methods using {M}-symmetry.
\newblock {\em Neural Parallel and Scientific Computations}, 20(3):437, 2012.

\bibitem{gottlieb2009}
Sigal Gottlieb, David~I. Ketcheson, and Chi-Wang Shu.
\newblock High order strong stability preserving time discretizations.
\newblock {\em Journal of Scientific Computing}, 38(3):251--289, 2009.

\bibitem{guibert2009cyclic}
D~Guibert and D~Tromeur-Dervout.
\newblock Cyclic distribution of pipelined parallel deferred correction method
  for {ODE/DAE}.
\newblock In {\em Parallel Computational Fluid Dynamics 2007}, pages 171--178.
  Springer, 2009.

\bibitem{Gustafsson_Kress_2001}
Bertil Gustafsson and Wendy Kress.
\newblock Deferred correction methods for initial value problems.
\newblock {\em BIT Numerical Mathematics}, 41(5):986–995, 2001.

\bibitem{hairer1993}
Ernst Hairer, Syvert~P. N\o~rsett, and G.~Wanner.
\newblock {\em Solving ordinary differential equations {I}: Nonstiff Problems}.
\newblock Springer Series in Computational Mathematics. Springer, Berlin,
  second edition, 1993.

\bibitem{Hosea1994a}
M.E. Hosea and L.F. Shampine.
\newblock {Efficiency comparisons of methods for integrating ODEs}.
\newblock {\em Computers \& Mathematics with Applications}, 28(6):45--55,
  September 1994.

\bibitem{Hull1972}
T.~E. Hull, W.~H. Enright, B.~M. Fellen, and A.~E. Sedgwick.
\newblock {Comparing Numerical Methods for Ordinary Differential Equations}.
\newblock {\em SIAM Journal on Numerical Analysis}, 9(4):603--637, December
  1972.

\bibitem{Ito_Fukushima_1997}
Takashi Ito and Toshio Fukushima.
\newblock Parallelized extrapolation method and its application to the orbital
  dynamics.
\newblock {\em The Astronomical Journal}, 114:1260, Sep 1997.

\bibitem{jackson}
K.~R. Jackson and S.~P.~N\o rsett.
\newblock The potential for parallelism in {R}unge-{K}utta methods . part 1 :
  Rk formulas in standard form.
\newblock 32(1):49--82, 2012.

\bibitem{Kappeller_Kiehl_Perzl_Lenke_1996}
M.~Kappeller, M.~Kiehl, M.~Perzl, and M.~Lenke.
\newblock Optimized extrapolation methods for parallel solution of {IVP}s on
  different computer architectures.
\newblock {\em Applied Mathematics and Computation}, 77(2-3):301–315, Jul
  1996.

\bibitem{Kennedy2000}
C.~A. Kennedy, M.~H. Carpenter, and R~Michael Lewis.
\newblock Low-storage, explicit {R}unge–-{K}utta schemes for the compressible
  {N}avier–-{S}tokes equations.
\newblock {\em Applied Numerical Mathematics}, 35(3):177--219, November 2000.

\bibitem{ketchesonleveque_periodic}
David~I. Ketcheson and Randall~J. LeVeque.
\newblock Shock dynamics in layered periodic media.
\newblock {\em Communications in Mathematical Sciences}, 10(3):859--874, 2012.

\bibitem{klp_internal}
David~I Ketcheson, Lajos L\'{o}czi, and Matteo Parsani.
\newblock Internal error propagation in explicit {R}unge--{K}utta
  discretization of {PDE}s.
\newblock Preprint available from \url{http://arxiv.org/abs/1309.1317}, 2013.

\bibitem{Ketcheson2011}
DI~Ketcheson, Matteo Parsani, and RJ~LeVeque.
\newblock {High-order Wave Propagation Algorithms for Hyperbolic Systems}.
\newblock {\em SIAM Journal on Scientific Computing}, 35(1):A351--A377, 2013.

\bibitem{liu2008}
Y~Liu, CW~Shu, and M~Zhang.
\newblock {Strong stability preserving property of the deferred correction time
  discretization}.
\newblock {\em Journal of Computational Mathematics}, 2007.

\bibitem{lustman1992solution}
L~Lustman, B~Neta, and W~Gragg.
\newblock Solution of ordinary differential initial value problems on an intel
  hypercube.
\newblock {\em Computers \& Mathematics with Applications}, 23(10):65--72,
  1992.

\bibitem{minion2010hybrid}
Michael Minion.
\newblock A hybrid parareal spectral deferred corrections method.
\newblock {\em Communications in Applied Mathematics and Computational
  Science}, 5:265--301, 2010.

\bibitem{ono200625}
Hiroshi Ono.
\newblock On the 25 stage 12th order explicit {R}unge-{K}utta method.
\newblock {\em Transactions-Japan Society for Industrial and Applied
  Mathematics}, 16(3):177, 2006.

\bibitem{Prince1981}
P.J. Prince and J.R. Dormand.
\newblock {High order embedded Runge-Kutta formulae}.
\newblock {\em Journal of Computational and Applied Mathematics}, 7(1):67--75,
  March 1981.

\bibitem{rauber1997load}
Thomas Rauber and Gudula R{\"u}nger.
\newblock Load balancing schemes for extrapolation methods.
\newblock {\em Concurrency: Practice and Experience}, 9(3):181--202, 1997.

\bibitem{Shampine1986}
L.~F. Shampine and L.~S. Baca.
\newblock {Fixed versus variable order Runge-Kutta}.
\newblock {\em ACM Transactions on Mathematical Software}, 12(1):1--23, May
  1986.

\bibitem{simonsen1990extrapolation}
HH~Simonsen.
\newblock {\em Extrapolation methods for ODE's: continuous approximations, a
  parallel approach}.
\newblock PhD thesis, PhD thesis, University of Trondheim, Norway, 1990.

\bibitem{VanderHouwen1990}
P.~J. van~der Houwen and B.~P. Sommeijer.
\newblock {Parallel ODE solvers}.
\newblock {\em ACM SIGARCH Computer Architecture News}, 18(3):71--81, September
  1990.

\bibitem{van1990parallel}
Peter~J Van Der~Houwen and Benjamin~P Sommeijer.
\newblock Parallel iteration of high-order runge-kutta methods with stepsize
  control.
\newblock {\em Journal of Computational and Applied Mathematics},
  29(1):111--127, 1990.

\end{thebibliography}

\end{document}